\documentclass[12pt]{amsart}
\usepackage[T1]{fontenc}
\usepackage[dvips]{graphicx}
\usepackage{amsmath,amsfonts,amssymb,amsthm,latexsym,graphicx}
\usepackage{a4wide,pdfsync}
\usepackage{hyperref,color}
\newtheorem{theo}{Theorem}
\newtheorem{lem}[theo]{Lemma}
\newtheorem{prop}[theo]{Proposition}
\newtheorem{defn}[theo]{Definition}

\newtheorem{rem}{Remark}

\newcommand\Z{\mathbb{Z}}
\newcommand\R{\mathbb{R}}
\newcommand\N{\mathbb{N}}

\newcommand\E{\mathbb{E}}
\def\d{\mathrm{\, d}}
\def\i{\mathbf{ i}}

\def\eqloi{\stackrel{\rm fidi}{=}} 

\author[Amblard]{Pierre-Olivier Amblard}
\address{Dept. Math \&
Stat.\\
The University of Melbourne\\
Parkville, VIC 3010, Australia\\
and\\
GIPSAlab/CNRS UMR 5216/ BP 46\\
38402 Saint Martin d'H\`eres cedex, France}
\email{pierre-olivier.amblard@gipsa-lab.grenoble-inp.fr}
\urladdr{http://www.gipsa-lab.grenoble-inp.fr/~bidou.amblard/}

\author[Coeurjolly]{Jean-Fran\c{c}ois Coeurjolly}
\address{Laboratoire Jean Kuntzmann, UMR 5224 \\
BSHM, 1251 Av. Centrale BP 47, \\ 
38040 Grenoble Cedex 9 FRANCE \\
and \\
GIPSAlab/CNRS UMR 5216/ BP 46\\ 
38402 Saint Martin d'H\`eres cedex, France}
\email{jean-francois.coeurjolly@upmf-grenoble.fr}
\urladdr{http://www-ljk.imag.fr/membres/Jean-Francois.Coeurjolly/}

\author[Lavancier]{Fr\'ed\'eric Lavancier }
\address{Universit\'e de Nantes\\ 
Laboratoire de math\'ematiques Jean
Leray\\
2, rue de la Houssini\`ere 44322 Nantes, France}
\email{frederic.lavancier@univ-nantes.fr}
\urladdr{http://www.math.sciences.univ-nantes.fr/~lavancie}

\author[Philippe]{Anne Philippe}
\address{Universit\'e de Nantes\\
Laboratoire de math\'ematiques Jean
Leray\\
2, rue de la Houssini\`ere 44322 Nantes, France}
\email{anne.philippe@univ-nantes.fr}
\urladdr{http://www.math.sciences.univ-nantes.fr/~philippe}

\title{Basic properties of the Multivariate  Fractional Brownian Motion}
\begin{document}
\begin{abstract}
This paper reviews and extends some recent results on the multivariate
fractional Brownian motion (mfBm) and its increment process.
A characterization of the mfBm through its covariance function is
obtained.  Similarly, the correlation and spectral analyses of the
increments are investigated.  On the other hand we show that (almost)
all   mfBm's  may be reached as the limit of partial sums of (super)linear processes. 
Finally, an algorithm to perfectly simulate the mfBm
is presented and illustrated by some simulations.
\end{abstract}

\thanks{Research supported  in part
by ANR STARAC grant and a fellowship from R\'egion Rh\^one-Alpes (France) and a Marie-Curie
International Outgoing Fellowship from the European Community.\\
Research supported in part by ANR InfoNetComaBrain grant.}
\maketitle

\section{Introduction}

\label{sec:introduction}
The fractional Brownian motion is the unique Gaussian self-similar process with
stationary increments. In the seminal paper of Mandelbrot and Van Ness
\cite{MandVN68}, many properties of the fBm and its increments are developed (see also \cite{Taq} for a review
of the basic properties). Depending on the scaling factor (called
Hurst parameter), the increment process may exhibit long-range dependence, and is commonly
used in modeling physical phenomena. However in many fields of applications (e.g. neuroscience, economy, sociology, physics, etc), multivariate
measurements are performed and they involve specific properties such as fractality, long-range dependence,
self-similarity, etc.  Examples can be found in economic time series
(see \cite{DavidHash}, \cite{FlemYHJ},
\cite{GilAlana}), genetic sequences \cite{AriCar09}, multipoint
velocity measurements in turbulence, functional Magnetic Resonance
Imaging of several regions of the brain \cite{AchaBMB08}.

It seems therefore natural to extend the fBm to a multivariate framework.
Recently, this question has been investigated in
\cite{LPS2009,LPS2010,CAA10}. The aim of this paper is to summarize
and to complete some of these advances on the multivariate fractional Brownian motion and its increments.
A multivariate extension of the fractional Brownian motion can be
stated as follows :
\begin{defn} \label{def-mfBm}
A Multivariate fractional Brownian motion ($p$-mfBm or mfBm) with parameter $H\in
(0,1)^p$ is a $p$-multivariate process
satisfying the  three following properties
 \begin{itemize}
 \item \label{item:1} Gaussianity,
 \item \label{item:2} Self-similarity with parameter $H\in(0,1)^p$,
 \item \label{item:3} Stationarity of the increments.
 \end{itemize}
\end{defn}
Here, self-similarity has to be understood as joint self-similarity. More formally, we use the following definition.
\begin{defn} \label{def-Hss}
 A multivariate process $(X(t) = (X_1(t), \cdots, X_p(t)))_{ t\in \R}
 $ is said self-similar if there exists a vector $H=(H_1, \cdots,
 H_p)\in (0,1)^p$ such that for any $\lambda >0$,
 \begin{equation} \label{def:selfsim} (X_1(\lambda t), \cdots,
  X_p(\lambda t)))_{t\in\R} \eqloi \big(\lambda^{H_1} X_1(t),
  \cdots, \lambda^{H_p} X_p(t)\big) _{t\in\R},
 \end{equation}
 where $\eqloi$ denotes the equality of finite-dimensional
 distributions. The parameter $H$ is called the self-similarity parameter.
\end{defn}
This definition can be viewed as a particular case of operator self-similar
processes by taking diagonal operators (see
\cite{DidierPipiras,HudsonMason,LahaRoh,MaejimaMAson}).

 Note that, as in the univariate case \cite{Lamperti62}, the Lamperti transformation induces an
 isometry between the self-similar and the stationary multivariate
 processes. Indeed, from Definition \ref{def-Hss}, it is not
 difficult to check that $(Y(t)) _{t\in\R} $ is a $p$-multivariate
 stationary process if and only if there exists $H\in(0,1)^p$ such
 that its Lamperti transformation $ (t^{H_1} Y_1(\log(t)) ,\ldots,t^{H_p} Y_p(\log(t)))
 _{t\in\R} $ is a $p$-multivariate self-similar process.

The paper is organized as follows.
 In Section \ref{sec:depend-struct-mfbm}, we study the
covariance structure of the mfBm and its increments. The
cross-covariance and the cross-spectral density of the increments lead
to interesting long-memory type properties. Section
\ref{sec:integr-repr} contains the time domain as well as
the spectral domain stochastic integral representations of the mfBm.
Thanks to these results we obtain a characterization of the mfBm
through its covariance matrix function. Section
\ref{limitprocess:sec} is devoted to limit theorems, the mfBm is
obtained as the limit of partial sums of linear processes. Finally, we discuss in Section \ref{simulation:sec} the problem of simulating sample paths of the mfBm.
We propose to use the Wood and Chan's algorithm \cite{WoodC94} well adapted to
generate multivariate stationary Gaussian random fields with
prescribed covariance matrix function.

\section{Dependence structure of the mfBm and of its increments }
\label{sec:depend-struct-mfbm}

\subsection{Covariance function of the mfBm}
\label{sec:covar-funct-mfbm}

In this part, we present the form of the covariance matrix of the mfBm.

 Firstly, as each component is a fractional brownian motion, the covariance function of the $i$-th component is well-known and we have
\begin{equation}
  \E X_i(s)X_i(t) \ = \frac{\sigma_i^2}{2} \left\{|s|^{2H_{i}} +
   |t|^{2H_{i}} - |t-s|^{2H_{i}}\right\}. \label{fbm}
 \end{equation}
with $\sigma^2_i := {\rm var}(X_i(1)) $. The cross covariances are given in the following proposition.

\begin{prop} [Lavancier \textit{et al.} \cite{LPS2009}] \label{prop-cross-cov}
 The cross covariances
of the mfBm satisfy the following representation, for all $(i, j)\in
\{1,\ldots,p\}^2 $, $ i\ne j$,
\begin{enumerate}
\item If $H_i + H_j \neq 1$, there exists $(\rho_{i,j},\eta_{i,j}
 )\in [-1,1] \times \R $ with $\rho_{i,j} = \rho_{j,i} =\mathrm{corr}(X_i(1),X_j(1))$ and $\eta_{i,j} = -\eta_{j,i}$ such that
\begin{multline} \E X_i(s) X_j(t) \ =
  \frac{\sigma_i\sigma_j}{2} \left\{ (\rho_{i,j}+\eta_{i,j}
    \mathrm{sign}(s) )  |s|^{H_i+ H_j} +
     (\rho_{i,j} - \eta_{i,j} \mathrm{sign}(t) )    |t|^{H_i+H_j}
    \right. \\
\left.     - ( \rho_{i,j} -\eta_{i,j} \mathrm{sign}(t-s) )  |t-s|^{H_i +
    H_j} \right\}.
\label{ccov-not1}
 \end{multline}
\item If $H_i+H_j = 1$, there exists $(\tilde\rho_{i,j}
 ,\tilde\eta_{i,j} ) \in [-1,1] \times \R $ with $\tilde\rho_{i,j} =
 \tilde\rho_{j,i} =\mathrm{corr}(X_i(1),X_j(1))$  and $\tilde\eta_{i,j} =
 -\tilde\eta_{j,i}$ such that
  \begin{equation}\label{ccov-1}
   \E X_i(s)X_j(t) \ = \frac{\sigma_i\sigma_j}{2} \left\{\tilde \rho_{i,j} (|s|+|t|-|s-t|)+ \tilde\eta_{i,j} (t\log |t|-s\log |s| - (t-s) \log |t-s|)\right \}.
  \end{equation}
\end{enumerate}
\end{prop}

\begin{proof}Under some conditions of regularity, Lavancier \textit{et
 al.} \cite{LPS2009} actually prove that Proposition 3 is true for any $L^2$ self-similar multivariate process with stationary increments. The form of cross covariances is obtained as the unique solution of a functional equation.
Formulae \eqref{ccov-not1} and \eqref{ccov-1} correspond to expressions
given in \cite{LPS2009} after the following reparameterization : $ \rho_{i,j} = (c_{i,j} +c_{j,i})/2 $ and $\eta_{i,j} = (c_{i,j} -
c_{j,i})/2$ where $c_{i,j}$ and $c_{j,i}$ arise in \cite{LPS2009}.

\end{proof}

\begin{rem}
Extending the definition of parameters
$\rho_{i,j},\tilde\rho_{i,j},\eta_{i,j},\tilde\eta_{i,j}$ to the case
$i=j$, we have $\rho_{i,i}
 =\tilde\rho_{i,i}=1$ and $\eta_{i,i} = \tilde\eta_{i,i} =0$, so that
  \eqref{fbm} coincides with
 \eqref{ccov-not1} and \eqref{ccov-1}.
\end{rem}

\begin{rem} \label{rem-constr}
The constraints on coefficients
$\rho_{i,j},\tilde\rho_{i,j},\eta_{i,j},\tilde\eta_{i,j}$ are
necessary but not sufficient conditions to ensure that the functions
defined by \eqref{ccov-not1} and \eqref{ccov-1} are covariance functions. This problem will be
discussed in Section \ref{sec:existence}.
\end{rem}
\begin{rem}\label{rem-exist}
Note that coefficients $\rho_{i,j},\tilde\rho_{i,j},\eta_{i,j},\tilde\eta_{i,j}$ depend on
the parameters $(H_i,H_j)$. Assuming the continuity of the cross
covariances function with respect to the parameters $(H_i,H_j)$, the
expression \eqref{ccov-1} can be deduced from \eqref{ccov-not1} by taking the limit as $H_i+H_j$
tends to $1$, noting that $((s+1)^H-s^H-1)/(1-H)\to s \log|s|-(s+1)\log|s+1|$ as $H\to 1$. We obtain the following relations between the
coefficients : as $ H_i+H_j\to 1$
$$\rho_{i,j} \sim \tilde\rho_{i,j} \qquad \text{ and }
\qquad (1-H_i-H_j)\eta_{i,j} \sim
\tilde\eta_{i,j}.$$
 This convergence result can suggest a reparameterization of
coefficients $\eta_{i,j} $ in $(1-H_i-H_j)\eta_{i,j}$.
\end{rem}

\subsection{ The increments process}

This part aims at exploring the covariance structure of the increments
of size $\delta$ of a multivariate fractional Brownian motion given by
Definition~\ref{def-mfBm}. Let $\Delta_\delta
X=(X(t+\delta)-X(t))_{t\in\R}$ denotes the increment process of the multivariate fractional Brownian motion of size $\delta$ and let $\Delta_\delta X_i$ be its $i$-th component.

Let $\gamma_{i,j}(h,\delta)=\E \Delta_\delta X_i(t) \Delta_\delta X_j(t+h) $ denotes the cross-covariance of the increments of size $\delta$ of the components $i$ and $j$. Let us introduce the function $w_{i,j}(h)$ given by
\begin{equation}
w_{i,j}(h) = \left\{ \begin{array}{ll}
(\rho_{i,j} -\eta_{i,j} \mathrm{sign}(h) ) |h|^{H_i+H_j} & \mbox{ if } H_i+H_j \neq 1, \\
\tilde{\rho}_{i,j} |h| + \tilde{\eta}_{i,j} h \log|h| & \mbox{ if } H_i+H_j =1.
\end{array} \right.
\end{equation}
Then from Proposition~\ref{prop-cross-cov}, we deduce that $\gamma_{i,j}(h,\delta)$ is given by
\begin{equation}\label{def-gi,j}
\gamma_{i,j}(h,\delta) = \frac{\sigma_i \sigma_j}2 \bigg( w_{i,j}(h-\delta)- 2 w_{i,j}(h) + w_{i,j}(h+\delta)\bigg).
\end{equation}

Now, let us present the asymptotic behaviour of the cross-covariance function.

\begin{prop} \label{prop-asymp-cross-cov}
As $|h|\to +\infty$, we have for any $\delta>0$
\begin{equation}
\gamma_{i,j}(h,\delta) \sim \sigma_i \sigma_j \delta^2 |h|^{H_i+H_j-2} \kappa_{i,j}(\mathrm{sign}(h)),
\end{equation}
with
\begin{equation}
\kappa_{i,j}(\mathrm{sign}(h))=\left\{
\begin{array}{ll}
(\rho_{i,j}-\eta_{i,j} \mathrm{sign}(h)) (H_i+H_j)(H_i+H_j-1) & \mbox{ if } H_i+H_j\neq 1, \\ 
\tilde\eta_{i,j} \mathrm{sign}(h)& \mbox{ if }H_i+H_j=1. 
\end{array}
\right.
\end{equation}
\end{prop}
\begin{proof}
Let $\alpha=H_i+H_j$. Let us choose $h$, such that $|h|\geq \delta$, which ensures that $\mathrm{sign}(h-\delta)=\mathrm{sign}(h)=\mathrm{sign}(h+\delta)$. When $\alpha\neq 1$, this allows us to write
$$
\gamma_{i,j}(h,\delta) = \frac{\sigma_i \sigma_j}2 |h|^\alpha \left(\rho_{i,j} -\eta_{i,j}\mathrm{sign}(h)\right) B(h),
$$
with $B(h)= \left(1-\frac\delta h\right)^\alpha-2+\left(1+\frac\delta h\right)^\alpha \sim \alpha(\alpha-1)\delta^2 h^{-2}$, as $|h| \to +\infty$. When $\alpha= 1$ and $|h|\geq \delta$, $\gamma_{i,j}(h,\delta)$ reduces to
$$
\gamma_{i,j}(h,\delta) = \frac{\sigma_i \sigma_j}2 \tilde{\eta}_{i,j} B(h) \mbox{ with } \;\;B(h)=\left( (h-\delta)\log\left(1-\frac\delta h\right) + (h+\delta)\log\left(1+\frac\delta h\right)\right).
$$
Using the expansion of $\log(1\pm x)$ as $x\to 0$ leads to $B(h)\sim \delta^2 |h|^{-1}$ as $|h|\to+\infty$, which implies the result.
\end{proof}

Proposition~\ref{prop-asymp-cross-cov} and \eqref{def-gi,j} lead to the following
important remarks on the dependence structure. For $i\neq j$ and $H_i+H_j \neq 1$ :
\begin{itemize}
\item If the two fractional Gaussian noises
 are short-range dependent (i.e. $H_i <1/2$ and $H_j<1/2$) then they are either short-range interdependent
 if $\rho_{i,j}\neq 0 $ or $ \eta_{i,j}\neq 0$, or independent if $\rho_{i,j}=\eta_{i,j}=0$.
\item If the two fractional Gaussian noises are long-range
 dependent (i.e. $H_i> 1/2$ and $H_j>1/2$) then they are either long-range
 interdependent if $\rho_{i,j}\neq 0 $ or $\eta_{i,j}\neq 0$, or independent if
 $\rho_{i,j}=\eta_{i,j}=0$. This confirms the dichotomy principle
 observed in \cite{DidierPipiras}.
\item In the other cases, the two fractional Gaussian
 noises can be short-range interdependent if $\rho_{i,j}\neq 0 $ or
 $\eta_{i,j}\neq 0$ and $H_i+H_j<1$, long-range interdependent if
 $\rho_{i,j}\neq 0 $ or $\eta_{i,j}\neq 0$ and $H_i+H_j>1$ or independent if $\rho_{i,j}=\eta_{i,j}=0$.
\end{itemize}
Moreover, note that when $H_i+H_j=1$, whatever the nature of the two fractional Gaussian noises (i.e. short-range or long-range dependent, or even independent), they are either long-range interdependent if $\tilde\eta_{i,j}\neq 0$ or independent if $\tilde\eta_{i,j}=0$.

The following result characterizes the spectral nature of the increments of a mfBm.

\begin{prop} [Coeurjolly \textit{et al.} \cite{CAA10}]\label{prop-spd}
 Let $S_{i,j}(\cdot ,\delta)$ be the (cross)-spectral density  of the increments of size $\delta$ of the components $i$ and $j$, \textit{i.e.} the Fourier transform of $\gamma_{i,j}(\cdot ,\delta)$
$$ S_{i,j}(\omega ,\delta)  = \frac{1}{2\pi} \int_\R e^{ - \i h \omega } \gamma_{i,j}(h ,\delta) \d h =: FT (\gamma_{i,j}(\cdot ,\delta)). $$

\noindent $(i)$ For all $i,j$ and for all $H_i,H_j$, we have
\begin{equation}
S_{i,j}(\omega,\delta) = \frac{\sigma_i\sigma_j}{\pi} \Gamma(H_i+H_j+1) \frac{1-\cos(\omega\delta)}{|\omega|^{H_i+H_j+1}} \times {\tau}_{i,j}(\mathrm{sign}(\omega)),
\label{eq:spectre}
\end{equation}
where
\begin{equation} \label{def-tautilde}
\tau_{i,j}(\mathrm{sign}(\omega)) = \left\{ \begin{array}{ll}
\rho_{i,j} \sin\left( \frac\pi 2(H_i+H_j)\right) - \i \eta_{i,j}
\mathrm{sign}(\omega)\cos\left( \frac\pi 2(H_i+H_j)\right) & \mbox{ if
} H_i+H_j \neq 1, \\ 
\tilde\rho_{i,j} - \i \frac{ \pi}2 \tilde\eta_{i,j}\mathrm{sign}(\omega)
& \mbox{ if } H_i+H_j = 1. 
\end{array} \right.
\end{equation}
$(ii)$ For any fixed $\delta$, when $H_i+H_j \neq 1$ then we have, as $\omega \to 0$,
\begin{equation}\label{eq-equiv:spd}
\big| S_{i,j}(\omega,\delta) \big| \sim \frac{\sigma_i\sigma_j}{2\pi} \Gamma(H_i+H_j+1)
   \delta^2 \; \frac{\left( \rho_{i,j}^2 \sin\left( \frac\pi 2(H_i+H_j)\right) ^2 + \eta_{i,j}^2
\cos\left( \frac\pi 2(H_i+H_j)\right) ^2\right)^{1/2}}{|\omega|^{H_i+H_j-1}}.
\end{equation}

$(iii)$ Moreover, when $H_i+H_j \neq 1$, the coherence function
between the two components $i$ and $j$ satisfies, for all $\omega$\\
\begin{eqnarray}
C_{i,j}(\omega,\delta) &: =& \frac{\left| S_{i,j}(\omega,\delta) \right|^2}{S_{i,i}(\omega,\delta) S_{j,j}(\omega,\delta) } \nonumber \\
&=&
 \frac{\Gamma(H_i+H_j+1)^2}{\Gamma(2H_i+1)\Gamma(2H_j+1)} \; \frac{
 \rho_{i,j}^2 \sin\left( \frac\pi 2(H_i+H_j)\right) ^2 + \eta_{i,j}^2
\cos\left( \frac\pi 2(H_i+H_j)\right) ^2 }{\sin(\pi H_i) \sin(\pi H_j)}.
  \label{eq-coherence:spd}
\end{eqnarray}

\noindent (iv) When $H_i+H_j=1$, \eqref{eq-equiv:spd} and \eqref{eq-coherence:spd} hold, replacing $ \rho_{i,j}^2 \sin\left( \frac\pi 2(H_i+H_j)\right) ^2 + \eta_{i,j}^2
\cos\left( \frac\pi 2(H_i+H_j)\right) ^2$ by
$\tilde\rho_{i,j}^2 + \frac{\pi^2}{4} \tilde\eta_{i,j}^2$.

\end{prop}

\begin{proof}
 The proof is essentially based on the fact that in the generalized
 function sense, for $\alpha>-1$,
 \begin{eqnarray*}
  FT(|h|^\alpha)&=& -\frac{1}{\pi}\Gamma(\alpha+1) \sin\left(\frac\pi 2 \alpha \right) |\omega|^{-\alpha-1}, \\
  FT(h_+^\alpha)&=& \frac{1}{2\pi} \Gamma(\alpha+1) e^{-\i
   \mathrm{sign}(\omega) \frac\pi 2 (\alpha+1)} |\omega|^{-\alpha-1} ,\\
  FT(h_-^\alpha)&=& \frac{1}{2\pi}\Gamma(\alpha+1) e^{\i \mathrm{sign}(\omega) \frac\pi 2 (\alpha+1)} |\omega|^{-\alpha-1} ,\\
  FT(h\log|h|) &=& \i \frac{\mathrm{sign}(\omega)}{2\omega^2}.
 \end{eqnarray*}
 See \cite{CAA10} for more details.
\end{proof}
\begin{rem}
From this proposition, we retrieve the same properties of dependence
and interdependence of $X_i$ and $X_j$ as stated after Proposition \ref{prop-asymp-cross-cov}.
\end{rem}

\subsection{Time reversibility}
\label{sec:reversibility}

A stochastic process is said to be time reversible if $X(t)=X(-t)$ for all $t$. As shows in \cite{DidierPipiras}, this is equivalent for zero-mean multivariate Gaussian stationary processes to
$\E X_i(t)X_j(s) = \E X_i(s)X_j(t)$ for $s,t\in \R$ or that the cross covariance of the increments satisfies $\gamma_{i,j}(h,\delta)=\gamma_{i,j}(-h,\delta)$ for $h\in\R$. The following proposition characterizes this property.

\begin{prop} \label{prop-rev}
A mfBm is time reversible if and only if $\eta_{i,j}=0$ (or $\tilde\eta_{i,j}=0$) for all $i,j=1,\ldots,p$.
\end{prop}

\begin{proof}
If $\eta_{i,j}=0$ (or $\tilde\eta_{i,j}=0$), $\gamma_{i,j}(h,\delta)$ is proportional to the covariance of a fractional Gaussian noise with Hurst parameter $(H_i+H_j)/2$ and is therefore symmetric. Let us prove the converse. Let $\alpha=H_i+H_j$, then
\begin{multline*}
\gamma_{i,j}(h,\delta)-\gamma_{i,j}(-h,\delta) = \sigma_i \sigma_j \times \\
\left\{ \begin{array}{ll}
- \eta_{i,j} \left( \mathrm{sign}(h-\delta)|h-\delta|^\alpha+2\mathrm{sign}(h)|h|^\alpha-\mathrm{sign}(h+\delta)|h+\delta|^\alpha\right) & \mbox{ if } \alpha\neq 1,\\
\tilde\eta_{i,j}\left(
(h-\delta)\log|h-\delta| -2h\log|h|+(h+\delta)\log|h+\delta| \right)& \mbox{ if } \alpha= 1.
\end{array}\right.
\end{multline*}
Assuming $\gamma_{i,j}(h,\delta)-\gamma_{i,j}(-h,\delta)$ equals zero
for all $h$ leads to $\eta_{i,j}=0$ (or $\tilde{\eta}_{i,j}=0$).
\end{proof}

\begin{rem}
This result can also be viewed from a spectral point view. The time
reversibility of a mfBm is equivalent to the fact that the spectral
density matrix is real. Using \eqref{eq:spectre}, this implies
$\eta_{i,j}=0$ (or $\tilde\eta_{i,j}=0$).
\end{rem}

\section{ Integral representation}
\label{sec:integr-repr}

\subsection{Spectral representation}
\label{sec:spectr-repr}
The following proposition contains the spectral representation of
mfBm. This representation will be especially
useful to obtain a condition easy to verify which ensures that the functions
defined by \eqref{ccov-not1} and \eqref{ccov-1} are covariance functions.
\begin{theo}[Didier and Pipiras, \cite{DidierPipiras}]
 \label{prop-spectraldomain}
Let $(X(t))_{t\in\R}$ be a mfBm with parameter $(H_1, \cdots,
H_p)\in (0,1)^p$. Then there exists a $p\times p$ complex matrix $A$ such
that  each component admits the following representation
\begin{equation}
 \label{eq:1}
 X_ i (t) = \sum_{j=1}^p \int \frac{e^{\i tx}-1}{\i x } (A_{ij} x_+^{-H_i +1/2}
 + \bar A_{ij} x_-^{-H_i +1/2} ) \tilde B_j(\d x) ,
\end{equation}
where for all $j=1,\dots,p$, $\tilde B_j$ is a Gaussian complex
measure such that $\tilde B_j= \tilde B_{j,1} + \i \tilde B_{j,2}$ with $\tilde
B_{j,1}(x)=\tilde B_{j,1}(-x)$, $\tilde B_{j,2}(x) =-\tilde B_{j,2}(x)$, $\tilde B_{j,1}$
and $\tilde B_{j,2}$ are independent and $E(\tilde B_{j,i}(\d x) \tilde B_{j,i}(\d x)
' ) = \d x$, $i=1,2$.

Conversely, any $p$-multivariate process satisfying \eqref{eq:1} is a
mfBm process.
\end{theo}
\begin{proof}
 This representation is deduced from the general spectral representation of operator fractional Brownian motions
obtained in \cite{DidierPipiras}. By denoting $
-\mathbb{H}+1/2 := {\rm diag}(-H_1+1/2, \cdots, -H_p +1/2)$
we have indeed
\begin{equation}\label{Xsp2}
X(t) = \int \frac{e^{\i tx}-1}{\i x } (x_+^{-\mathbb H +1/2 }A
 +  x_-^{-\mathbb H+1/2} \bar A) \tilde B(\d x).
\end{equation}

\end{proof}
Any mfBm having representation \eqref{eq:1} has a covariance function as in Proposition \ref{prop-cross-cov}. The
 coefficients $\rho_{i,j}$, $\eta_{i,j}$, $\tilde\rho_{i,j}$ and $\tilde\eta_{i,j}$ involved in \eqref{ccov-not1} and \eqref{ccov-1}
satisfy
\begin{equation}\label{lien-spect-cov}
(AA^*)_{i,j}= \frac{\sigma_i\sigma_j}{2\pi} \Gamma(H_i+H_j+1){\tau}_{i,j}(1),
\end{equation}
where $\tau_{i,j}$ is given in \eqref{def-tautilde} and where $A^*$
is the transpose matrix of $\bar A$. This relation is obtained by identification of the spectral matrix of the increments deduced on the one hand from \eqref{eq:1} and provided on the other hand in Proposition \ref{prop-spd}.

Given \eqref{eq:1}, relation \eqref{lien-spect-cov} provides easily
the coefficients $\rho_{i,j}$, $\eta_{i,j}$, $\tilde\rho_{i,j}$ and
$\tilde\eta_{i,j}$ which define the covariance function. The converse
is more difficult to obtain. Given a covariance function as in
Proposition \ref{prop-cross-cov}, obtaining the explicit
representation \eqref{eq:1} requires finding a matrix $A$ such that
\eqref{lien-spect-cov} holds. This choice is possible if and only if
the matrix on the right hand side of \eqref{lien-spect-cov} is
positive semidefinite. Then a matrix $A$ (which is not unique) may be deduced by the Cholesky decomposition. When $p=2$, an explicit solution is the matrix with entries, for $i,j=1,2$,
 \begin{multline*}A_{i,j}=\lambda_{i,j} \left[\left(\rho_{i,j}\sin\left(\frac{\pi}{2}(H_i+H_j)\right)+\eta_{i,j}\sqrt{\frac{1-C_{i,j}}{C_{i,j}}}\cos\left(\frac{\pi}{2}(H_i+H_j)\right)\right) \right. \\
\left.  + i \left(\rho_{i,j}\sqrt{\frac{1-C_{i,j}}{C_{i,j}}}\sin\left(\frac{\pi}{2}(H_i+H_j)\right)- \eta_{i,j}\cos\left(\frac{\pi}{2}(H_i+H_j)\right)\right)\right],\end{multline*}
where $\lambda_{i,j}=\dfrac{\sigma_i}{2\sqrt{\pi}}\dfrac{\Gamma(H_i+H_j+1)}{\sqrt{ \Gamma(2H_j+1)\sin(H_j\pi)}}$ and $C_{i,j}$ is given in \eqref{eq-coherence:spd}, provided $H_1+H_2\not=1$. When $H_1+H_2=1$, the same solution holds, replacing $\rho_{i,j}$ by $\tilde\rho_{i,j}$ and $\eta_{i,j}\cos\left(\frac{\pi}{2}(H_1+H_2)\right)$ by $-\frac{\pi}{2}\tilde\eta_{i,j}$.

\subsection{Moving average representation}
\label{sec:moving-aver-repr}

In the next proposition, we give an alternative characterization of
the mfBm from an integral representation in the time domain (or moving average representation).

\begin{theo}[Didier and Pipiras, \cite{DidierPipiras}]\label{prop-marepresentation}
 Let $(X(t))_{t\in\R}$ be a mfBm with parameter $(H_1, \cdots,
H_p)\in(0,1)^p$. Assume that for all $i\in\{1,...,p\}$,  $H_i\not =
1/2 $. Then there exist $M^+, M^-$ two $p\times p$ real matrices such
that  each component admits the following representation
\begin{equation}\label{Xma}
X_i(t)=\sum_{j=1}^p   \int_{\R} M_{i,j}^+ \left( (t-x)_+^{H_i-.5}-
  (-x)_+^{H_i-.5}\right) +
 M_{i,j}^- \left((t-x)_-^{H_i-.5} - (-x)_-^{H_i-.5}\right) W_j( \mathrm{d} x),
\end{equation}
with  $W(\mathrm{d}x) = (W_1(\mathrm{d}x), \cdots, W_p(\mathrm{d}x)) $ is a Gaussian white noise with zero mean, independent components and covariance
$\E W_i(\mathrm{d}x) W_j(\mathrm{d}x) = \delta_{i,j} \mathrm{d} x $.

Conversely, any $p$-multivariate process satisfying \eqref{Xma} is a
mfBm process.

\end{theo}

\begin{proof}
 This representation is deduced from the general representation
obtained in \cite{DidierPipiras}.
\end{proof}

\begin{rem} When  $H_i=1/2$ for each $i\in\{1,...,p\}$, it is shown
 in \cite{DidierPipiras} that  each component
of the mfBm admits the following representation :
 \begin{equation*}
  X_i(t)=\sum_{j=1}^p \int_{\R} M_{i,j}^+ (\mathrm{sign}(t-x)-\mathrm{sign}(x)) +
  M_{i,j}^- \left(\log|t-x| - \log|x|\right) W_j( \mathrm{d} x).
   \end{equation*}
Our conjecture is that this representation remains valid
when $H_i =1/2$ whatever the values of other parameters $H_j$, $j\not = i$.
\end{rem}


Starting from the moving average
representation \eqref{Xma}, using results in \cite{stoevTaqqu}, we can specify the coefficients $\rho_{i,j}$, $\eta_{i,j}$, $\tilde\rho_{i,j}$ and $\tilde\eta_{i,j}$ involved in the covariances \eqref{ccov-not1} and \eqref{ccov-1} (see \cite{LPS2009}). More precisely,
let us denote
$$
M^+ (M^{+}){'} = \left(\alpha^{++}_{i,j}\right), \quad M^- (M^{-}){'} =
\left(\alpha^{--}_{i,j}\right), \quad M^+ (M^{-}){'} =
\left(\alpha^{+-}_{i,j}\right)
$$
where $M'$ is the transpose matrix of $M$.
The variance of each component is equal to
 $$\sigma^2_i = \frac{B(H_i+.5, H_i +.5)}{\sin (H_i
    \pi)} \left\{ \alpha^{++}_{i,i} + \alpha^{--}_{i,i}
   - 2\sin (H_i\pi) \alpha^{+-}_{i,i} \right\},$$
where $B(\cdot , \cdot)$ denotes the Beta function. \\
Moreover, if ${H_i+H_j \not= 1}$ then
\begin{align*}
   {\sigma_i\sigma_j} \rho_{i,j}&=\frac{B(H_i+.5, H_j +.5)}{\sin ((H_i+H_j)
    \pi)} \times \\
   & \quad  \left\{ ( \alpha^{++}_{i,j}+\alpha^{--}_{i,j}) (\cos(H_i
    \pi) + \cos(H_j \pi) )  - (\alpha^{+-}_{i,j} + \alpha^{-+}_{i,j}) \sin ((H_i+H_j)\pi)\right\},\\
   {\sigma_i\sigma_j} \eta_{i,j}&=\frac{B(H_i+.5, H_j +.5)}{\sin ((H_i+H_j)
    \pi)} \times  \\
    &\quad \left\{ ( \alpha^{++}_{i,j} - \alpha^{--}_{i,j}) (\cos(H_i
    \pi) -  \cos(H_j \pi) )  - (\alpha^{+-}_{i,j} - \alpha^{-+}_{i,j}) \sin ((H_i+H_j)\pi)\right\}.
  \end{align*}

If ${H_i+H_j= 1}$ then
 \begin{align*}
   {\sigma_i\sigma_j} \tilde\rho_{i,j}&=B(H_i+.5, H_j+.5)
   \bigg\{\frac{\sin(H_i\pi)+\sin(H_j\pi)}{2}(\alpha^{++}_{i,j} + \alpha^{--}_{i,j}) - \alpha^{+-}_{i,j} - \alpha^{-+}_{i,j} \bigg\}, \\
    {\sigma_i\sigma_j} \tilde\eta_{i,j}&=
(H_j-H_i)(\alpha^{++}_{i,j} - \alpha^{--}_{i,j}).
  \end{align*}

Conversely, given a covariance function as in Proposition
\ref{prop-cross-cov}, if $H_i\not=1/2$ for all $i$, one may find
matrices $M^+$ and $M^-$ such that \eqref{Xma} holds, provided the
matrix on the right hand side of \eqref{lien-spect-cov} is positive semidefinite. Indeed, in this case, a matrix $A$ which solves
\eqref{lien-spect-cov} may be found by the Cholesky decomposition, then $M^+$ and $M^-$ are deduced from relation (3.20) in \cite{DidierPipiras}:

$$M^{\pm}=\sqrt{\frac{\pi}{2}}\left(D_1^{-1}A_1\pm D_2^{-1} A_2\right),$$
where $A=A_1+\i A_2$ and $$D_1=\mathrm{diag}\left(\sin(\pi H_1)\Gamma(H_1+\frac{1}{2}),\dots,\sin(\pi H_p)\Gamma(H_p+\frac{1}{2})\right),$$ $$D_2=\mathrm{diag}\left(\cos(\pi H_1)\Gamma(H_1+\frac{1}{2}),\dots,\cos(\pi H_p)\Gamma(H_p+\frac{1}{2})\right).$$

\subsection{Two particular examples} \label{sec-examples}

Let us focus on two particular examples which are quite natural: the causal mfBm ($M^-=0$) and the well-balanced mfBm ($M^-=M^+$). In the causal case, the integral representation is a direct generalization of the integral representation of Mandelbrot and Van Ness \cite{MandVN68} to the multivariate case. The well-balanced case is studied by Stoev and Taqqu in one dimension \cite{stoevTaqqu}. With the notation of the two previous sections, we note that the causal case (resp. well-balanced case) leads to $A_1 = \tan(\pi H) A_2$ (resp. $A_2=0$), where $\tan(\pi H):=\mathrm{diag}\left( \tan(\pi H_1),\ldots,\tan(\pi H_p)\right)$.
In these two cases, the covariance only depends on one parameter, for
instance $\rho_{i,j}$ (or $\tilde\rho_{i,j}$). Indeed we easily deduce $\eta_{i,j}$
(or $\tilde\eta_{i,j}$) as follows :
\begin{itemize}
\item in the causal case \textit{i.e.} $M^-=0$ or equivalently $
  A_1=\tan(\pi H) A_2$ :
\begin{eqnarray*}
\eta_{i,j} &=& - \rho_{i,j}
\;\tan(\frac\pi2(H_i+H_j))\tan(\frac\pi2(H_i-H_j)) \quad \text{ if }
H_i+H_j\neq 1,\\
\tilde{\eta}_{i,j} &=& \tilde{\rho}_{i,j} \frac{2}{\pi \tan(\pi H_i)}
\quad \text{ if }
H_i+H_j = 1.\\
\end{eqnarray*}
\item in the well-balanced case \textit{i.e.} $M^-=M^+$ or
  equivalently $A_2=0$ :
\begin{eqnarray*}
\eta_{i,j} &=& 0 \quad \text{ if }
H_i+H_j\neq 1, \\
\tilde\eta_{i,j} &=& 0 \quad \text{ if }
H_i+H_j = 1.\\
\end{eqnarray*}
\end{itemize}
\begin{rem}
 From Proposition~\ref{prop-rev}, the well-balanced
mfBm is therefore time reversible.
\end{rem}
\subsection{Existence of the covariance of the mfBm}
\label{sec:existence}

In this paragraph, we highlight some of the previous results in order
to exhibit the
sets of the possible parameters $(\rho_{i,j},\eta_{i,j})$ or
$(\tilde\rho_{i,j},\tilde\eta_{i,j}$) ensuring the existence of the
covariance of the mfBm.

For $i,j=1,\ldots,p$, let us give $(H_i,H_j)\in(0,1)^2$,
$(\sigma_i,\sigma_j)\in\R^+\times\R^+$ and $(\rho_{i,j}$
$\eta_{i,j}) \in [-1,1]\times \R$ with $\rho_{j,i}=\rho_{i,j}$ and $\eta_{j,i}=-\eta_{i,j}$ if
$H_i+H_j\neq 1$, or $(\tilde\rho_{i,j}$, $\tilde\eta_{i,j})\in
[-1,1]\times \R$ with $\tilde\rho_{j,i}=\tilde\rho_{i,j}$ and $\tilde\eta_{j,i}=-\tilde\eta_{i,j}$ if
$H_i+H_j=1$.

For this set of parameters, let us define the matrix $\Sigma(s,t)=(\Sigma_{i,j}(s,t))$ as follows :  $\Sigma_{i,i}(s,t)$ is given by \eqref{fbm} and
$\Sigma_{i,j}(s,t)$ is given by \eqref{ccov-not1} when $H_i+H_j \neq
1$ and
\eqref{ccov-1} when $H_i+H_j = 1$.
\begin{prop} \label{prop-exist} ${ }$
The matrix $\Sigma(s,t)$ is a covariance matrix function if and only
if the
Hermitian matrix $Q=\left( \Gamma(H_i+H_j+1){\tau}_{i,j}(1)
\right)$ with $\tau_{i,j}$ defined in
\eqref{def-tautilde}, is positive semidefinite. When $p=2$, this condition reduces to $C_{1,2} \leq 1$ where $C_{1,2}$ is the coherence defined by~(\ref{eq-coherence:spd}).
\end{prop} 

\begin{proof}
First, note that since $\rho_{j,i}=\rho_{i,j}$ and
$\eta_{j,i}=-\eta_{i,j}$ , $Q$ is a Hermitian matrix. Now, if $Q$ is
positive semidefinite, then so is the matrix $(2\pi)^{-1} (\sigma_i
\sigma_j Q_{i,j})$. Therefore there exists a matrix $A$
satisfying~(\ref{lien-spect-cov}). From Theorem
\ref{prop-spectraldomain}, there exists a mfBm having $\Sigma(s,t)$ as
covariance matrix function.
Conversely, if $\Sigma(s,t)$ is a covariance matrix function of a mfBm
then the representation \eqref{eq:1} holds and by
\eqref{lien-spect-cov}, the matrix $Q$ is positive semidefinite.

When $p=2$, the result comes from the fact that $Q$ is positive
semidefinite if and only if $\det(Q)\geq 0$ or equivalently $C_{1,2}
\leq 1$.
\end{proof}
\begin{figure}[h]
  \begin{center}
    \begin{tabular}{ll} {\includegraphics[scale=.45]{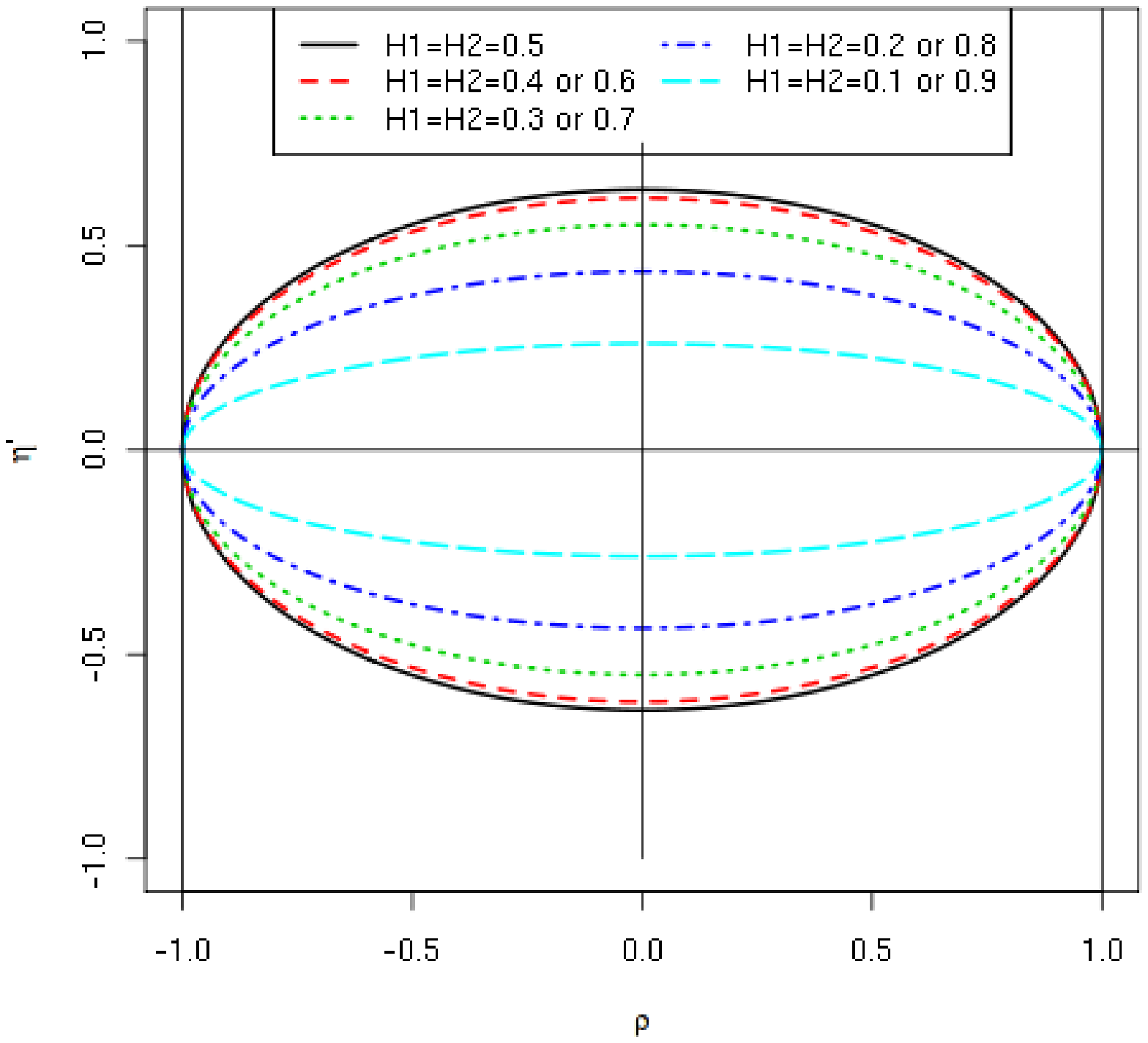}} &
      {\includegraphics[scale=.45]{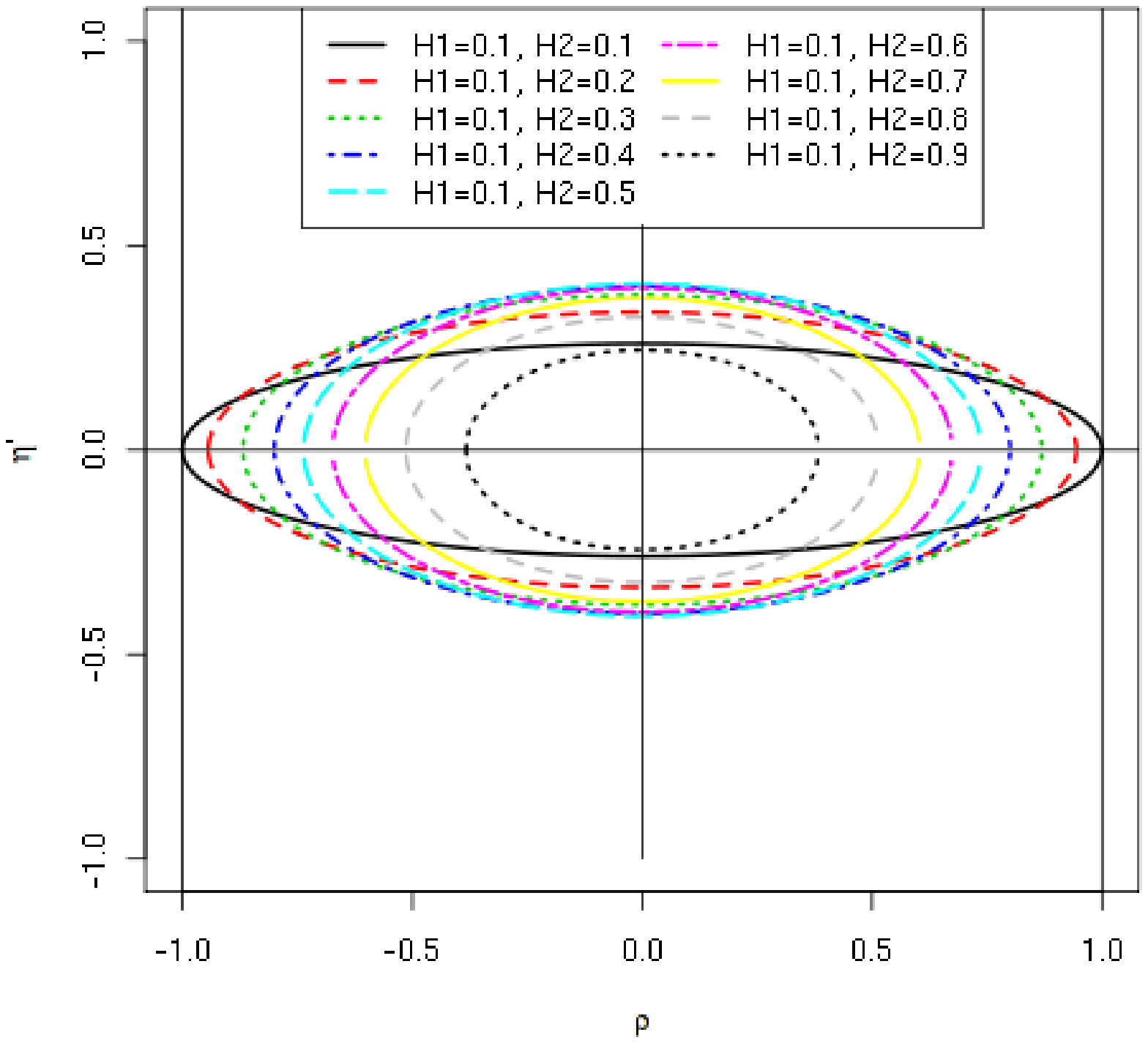}} \\
      {\includegraphics[scale=.45]{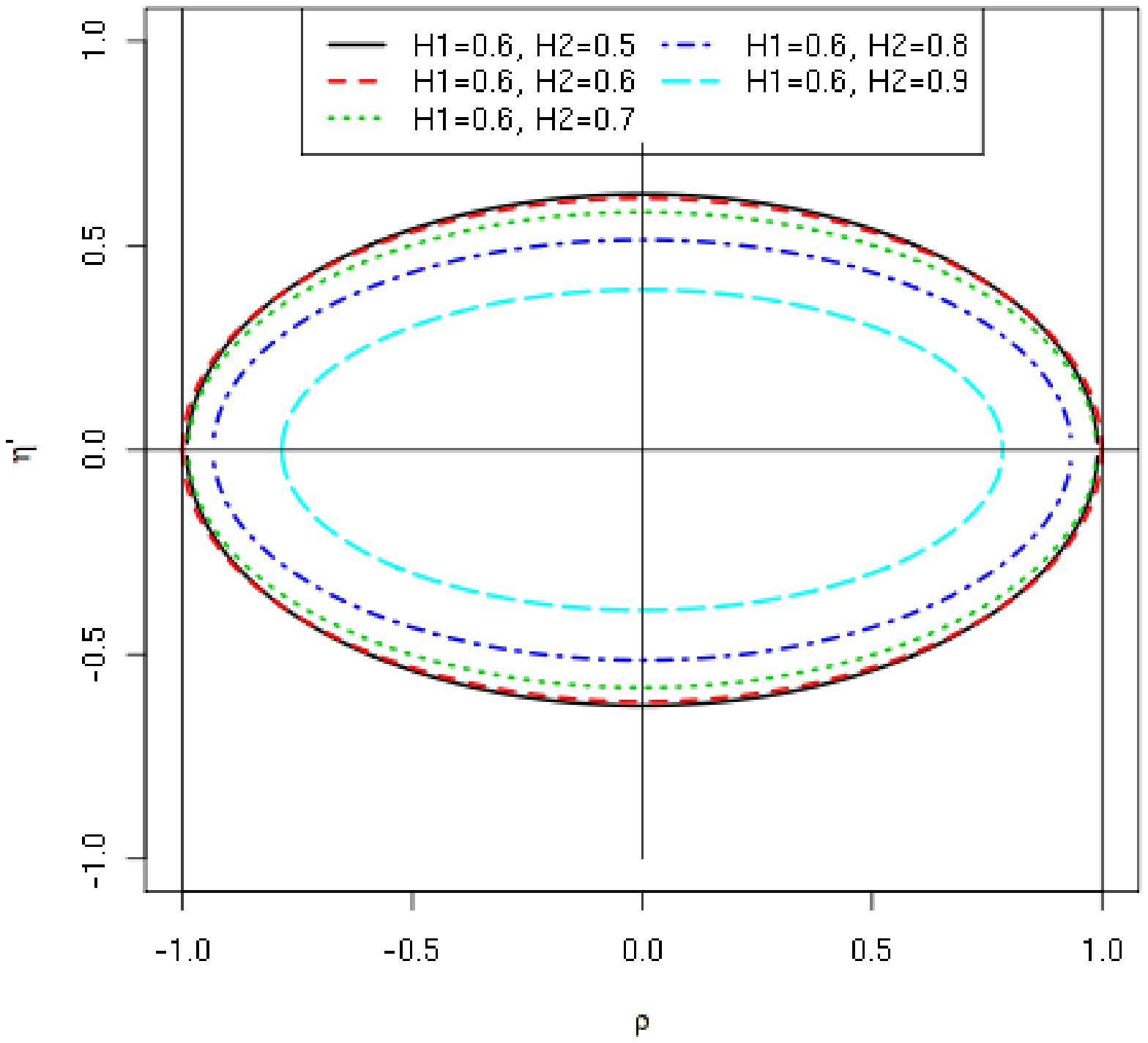}}
      & \includegraphics[scale=.45]{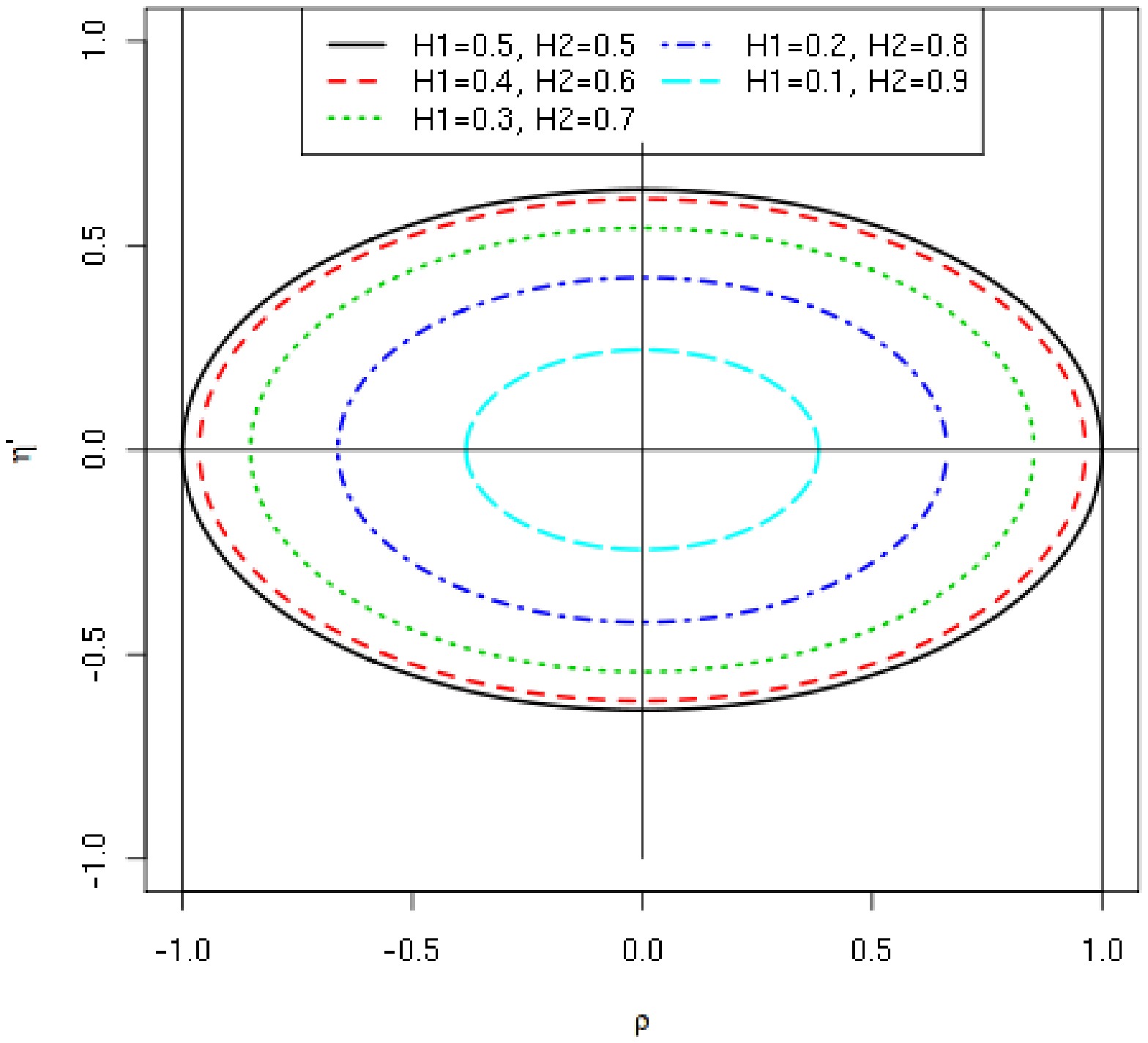}
    \end{tabular}
  \end{center}
\caption{Various examples of possible values for
 $(\rho_{1,2},\eta_{1,2}^\prime)$ with
{ $\eta_{1,2}^\prime:= \eta_{1,2} (1-H_1-H_2) $}  when $H_1+H_2\neq 1$ and
 $(\tilde\rho_{1,2},\tilde\eta_{1,2})$ when  {$H_1+H_2= 1$},
 ensuring that $\Sigma(s,t)$ is a covariance matrix function in the particular case $p=2$. }
\label{fig-ell}
\end{figure}

\begin{center}
\begin{figure}[h]
  \begin{center}
    \begin{tabular}{cc} {\includegraphics[scale=.45]{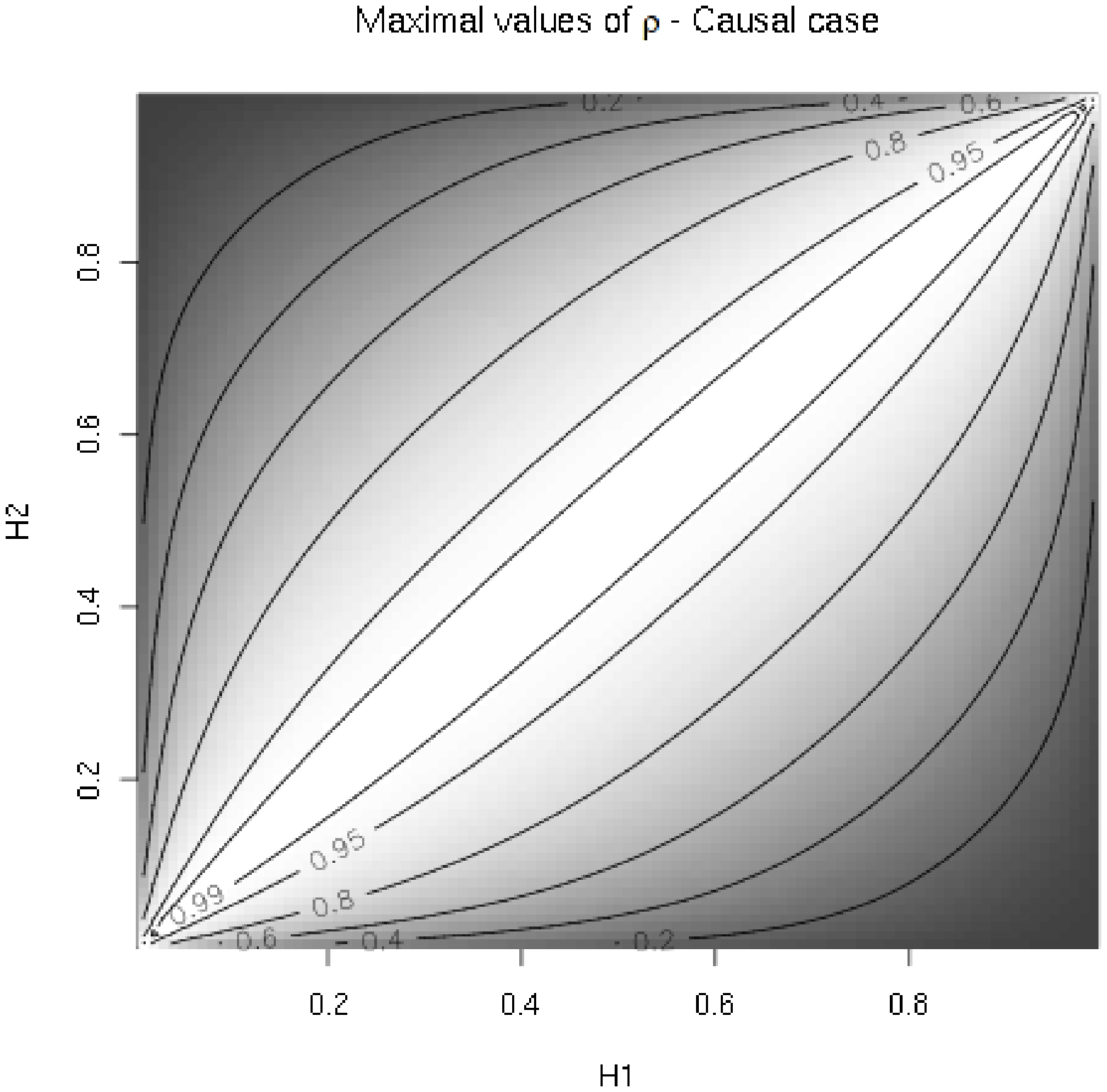}} &
      {\includegraphics[scale=.45]{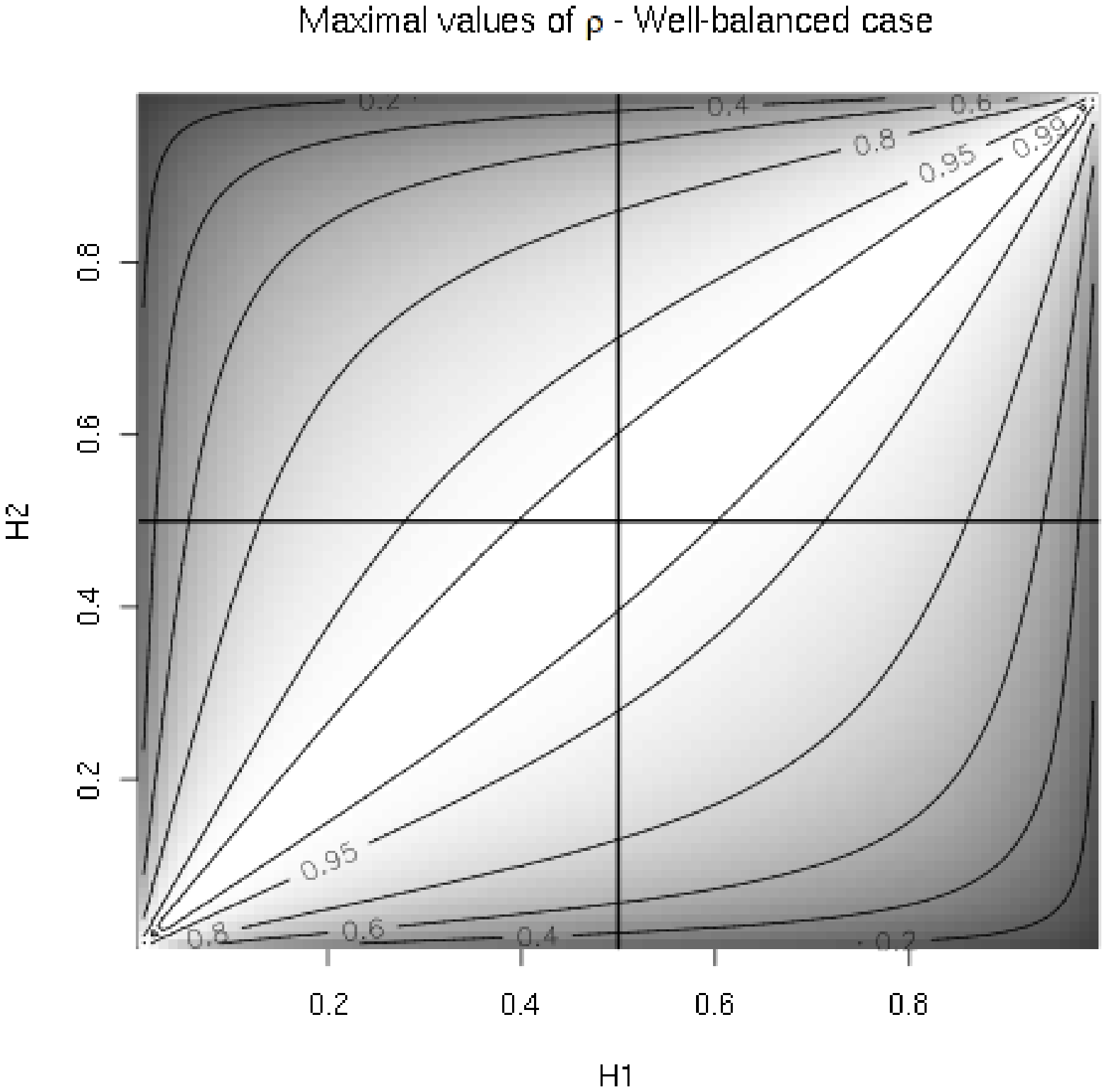}}
    \end{tabular}
  \end{center}

\caption{\small Maximal values of the absolute possible correlation
 parameter $|\rho_{1,2}|$ ensuring that $\Sigma(s,t)$ is a covariance
 matrix function in the case p = 2, in terms of $H_1$ and $H_2$ for the causal and well-balanced mfBm.}
\label{fig-detA}
\end{figure}
\end{center}

When $p=2$, for fixed values of $(H_1,H_2)$ the condition $C_{1,2}
\leq 1$ means that the set of possible parameters
$(\rho_{1,2},\eta_{1,2})$ is the interior of an ellipse. These sets
are represented in Figure \ref{fig-ell} according to different values
of $H_1$ and $H_2$. Note that, in order to compare
the cases $H_1+H_2\neq 1$ and $H_1+H_2=1$, we have reparameterized $\eta_{1,2}$
by {$\eta_{1,2}^\prime:=\eta_{1,2}(1-H_1-H_2)$}. In such a way, the second ellipse
becomes the limit of the first one as $H_1+H_2\to 1$ (see also Remark \ref{rem-exist}).

Let us underline that the maximum possible correlation between two fBm's is
obtained when $\eta_{1,2}=0$, i.e. when the $2$-mfBm is time reversible
according to Proposition \ref{prop-rev}.

\begin{rem}
When $H_1=\ldots=H_p=H\neq 1/2$, the matrix $Q$ rewrites
$Q_{i,j}=\Gamma(2H+1)( \sin(\pi H) \rho_{i,j} - \i \eta_{i,j}\cos(\pi H))$ and
\begin{itemize}
\item if the mfBm is time reversible, i.e. $\eta_{i,j}=0$ (for
 $i,j=1,\ldots,p$), then $Q$ is a correlation matrix and is therefore
 positive-semidefinite for any $|\rho_{i,j}|\leq 1$,
\item when $p=2$, the set of possible values for
 $(\rho_{1,2},\eta_{1,2})$ associated to $H$ and $1-H$ are the same.
\end{itemize}
\end{rem}

In the particular case of the causal or the well-balanced mfBm, the matrix
$\Sigma(s,t)$ can be expressed through the sole parameter
$\rho_{i,j}$. Let $\lambda(H_1,H_2)$ the function which equals to $\cos(\frac{\pi}2(H_1-H_2))^2$ in the causal case and which equals 1 in the well-balanced case. The maximal possible correlation when $p=2$ is given by
\begin{equation*} 
\rho_{1,2}^2 = \frac{\Gamma(2H_1+1)\Gamma(2H_2+1)}{\Gamma(H_1+H_2+1)^2} \frac{\sin(\pi H_1)\sin(\pi H_2)}{\sin(\frac\pi 2(H_1+H_2))^2}  \times \lambda(H_1,H_2).
\end{equation*}
Figure~\ref{fig-detA} represents $|\rho_{1,2}|$ with respect to $(H_1,H_2)$.

Figures~\ref{fig-ell} and~\ref{fig-detA} illustrate the main
limitation of the mfBm model. Under self-similarity condition
\eqref{def:selfsim}, it is not possible to construct arbitrary correlated
fractional Brownian motions. For example, when $H_1=0.1$
and $H_2=0.8$, the correlation cannot exceed $0.514$.

\section{The mfBm as a limiting process. }
\label{limitprocess:sec}

A natural way to construct self-similar processes is through limits of stochastic
 processes. In dimension one, the result is due to Lamperti \cite{Lamperti62}.
In \cite{HudsonMason}, an extension to operator self-similar
processes is given. A similar result for the mfBm is deduced and stated below.
In the following, a $p$-multivariate process $(X(t))_{t\in\R}$ is said proper if, for each $t$, the law of $X(t)$ is not contained in a proper subspace of $\R^p$.

\begin{theo}\label{lamperti}
Let $(X(t))_{t\in\R}$ be a $p$-multivariate proper process, continuous in
probability. If there exist a $p$-multivariate process $(Y(t))_{t\in\R}$ and
$p$ real functions $a_1,....,a_p$ such that
\begin{equation}
 \left(a_1(n) Y_1(nt), \ldots,a_p(n) Y_p(nt)\right)
\xrightarrow[fidi]{n\to \infty} X(t), \label{eq:3}
\end{equation}
then the multivariate process $(X(t))$ is self-similar.
Conversely, any multivariate self-similar process can be obtained as a
such limit.

\end{theo}

\begin{proof}

 The proof is similar to Theorem 5 in \cite{HudsonMason}.
Fix $k\in \N$ and $r>0$. For each $T\in\R^k$ we
denote $X(T): = (X(T_1) ,\ldots , X(T_k))$. Let $\mathcal{D}_{r,k} $ be
the set of all invertible diagonal matrices $\alpha$ such that, for all $T\in\R^k$,
$X(rT) =\alpha X(T)$.

Let us first show that $\mathcal{D}_{r,k} $ is not empty. According to
(\ref{eq:3}), we have
$$
 \mathrm{diag}(a_1(n) , \ldots,a_p(n) ) Y(nrT)
\xrightarrow[d]{n\to \infty} X(rT),
$$
and
$$
 \mathrm{diag}(a_1(rn) , \ldots,a_p(rn) ) Y(nrT)
\xrightarrow[d]{n\to \infty} X(T).
$$
Since $(X(t))$ is proper, $ \mathrm{diag}(a_1(n) , \ldots,a_p(n) ) $
and $ \mathrm{diag}(a_1(rn) , \ldots,a_p(rn) ) $
are invertible for $n$ large enough. Then, Theorem 2.3 in
\cite{Weissman} ensures that $\alpha_n$ defined by
$$\alpha_n = \mathrm{diag}(a_1(n) ,
\ldots,a_p(n) ) \mathrm{diag}(a_1(nr) ,
\ldots,a_p(nr) )^{-1} $$ has a limit in $\mathcal{D}_{r,k}
$. Moreover if $\alpha$ is a limit of $\alpha_n$ then $ X(rT) = \alpha
X(T) $ and thus $\mathcal{D}_{r,k} \not = \emptyset$.

It is then straightforward to adapt Lemma 7.2-7.5  in
\cite{HudsonMason} for the subgroup $\mathcal{D}_{r,k}$, which yields that for each $r$, $\cap_k \mathcal{D}_{r,k}$ is not empty.
Therefore, for any fixed $r>0$, there exists $\alpha \in\cap_k \mathcal{D}_{r,k} $
such that $ (X(rt))$ and $(\alpha X(T)) $ have the same finite
dimensional distributions. Theorem 1 in \cite{HudsonMason} ensures
that there exists $(H_1,\ldots,H_p)\in(0,1)^p$ such that $\alpha=\mathrm{diag}(r^{H_1},\ldots,r^{H_p})$.
The converse is trivial.

\end{proof}

As an illustration of Theorem \ref{lamperti}, the mfBm can be obtained as the weak limit of partial
sums of sum of linear processes (also called superlinear processes,
see \cite{zhao}). Some examples may be found in \cite{chung} and \cite{LPS2010}. In Proposition \ref{partialsums} below, we give a general convergence result which allows to reach almost any mfBm from such partial sums. The unique restriction concerns the particular case when at least one of the Hurst parameters is equal to $1/2$.

Let $(\epsilon_j(k))_{k\in Z}$, $j=1,\dots,p$ be $p$ independent
i.i.d. sequences with zero mean and unit variance. Let us consider the
superlinear processes
\begin{equation}
Z_i(t)=\sum_{j=1}^p \sum_{k\in\Z} \psi_{i,j}(t-k)\epsilon_j(k),\quad i=1,\dots,p,
\end{equation}
where $\psi_{i,j}(k)$ are real coefficients with $\sum_{k\in\Z} \psi_{i,j}^2(k)<\infty$.

Moreover, we assume that $\psi_{i,j}(k)=\psi_{i,j}^+ (k)+\psi_{i,j}^-(k)$ where $\psi_{i,j}^+(k)$ satisfies one of the following conditions:
\begin{itemize}
\item[(i)] $\psi_{i,j}^+(k)=\left(\alpha_{i,j}^+ + o(1)\right) k_+^{d_{i,j}^+ - 1}$ as $|k|\to\infty$, with $0<d_{i,j}^+<\frac{1}{2}$ and $\alpha_{i,j}^+\not=0$,
\item[(ii)] $\psi_{i,j}^+(k)=\left(\alpha_{i,j}^+ + o(1)\right) k_+^{d_{i,j}^+ - 1}$ as $|k|\to\infty$, with $-\frac{1}{2}<d_{i,j}^+<0$, $\sum_{k\in\Z} \psi_{i,j}^+(k)=0$ and $\alpha_{i,j}^+\not=0$,
\item[(iii)] $\sum_{k\in\Z} \left|\psi_{i,j}^+(k)\right|<\infty$ and let $\alpha_{i,j}^+:=\sum_{k\in\Z} \psi_{i,j}^+(k)\not=0$, $d_{i,j}^+:=0$.
\end{itemize}
Similarly, $\psi_{i,j}^-(k)$ is assumed to satisfy $(i)$, $(ii)$ or $(iii)$ where $k_+$, $d_{i,j}^+$ and $\alpha_{i,j}^+$ are replaced by $k_-$, $d_{i,j}^-$ and $\alpha_{i,j}^-$.

\begin{prop}\label{partialsums}
Let $d_i=max(d_{i1}^+,d_{i1}^-,\cdots,d_{ip}^+,d_{ip}^-)$, for $i=1,\dots,p$. Consider the vector of partial sums, for $\tau\in\R$,
$$S_n(\tau)=\left(n^{-d_1-(1/2)}\sum_{t=1}^{[n\tau]}Z_1(t),\cdots,n^{-d_p-(1/2)}\sum_{t=1}^{[n\tau]}Z_p(t)\right).$$
Then the finite dimensional distributions of $(S_n(\tau))_{\tau\in\R}$ converge in law towards a $p$-mfBm $(X(\tau))_{\tau\in\R}$.
\begin{itemize}
\item When $d_i\not=0$, $(X_i(\tau))_{\tau\in\R}$ is defined through the integral representation (\ref{Xma}) where $M_{i,j}^+=\alpha_{i,j}^+d_i^{-1}\mathbf{1}_{d_{i,j}^+=d_i}$ and $M_{i,j}^-=\alpha_{i,j}^-d_i^{-1}\mathbf{1}_{d_{i,j}^-=d_i}$.
\item When $d_i=0$, $X_i(\tau)=\sum_{j=1}^p (\alpha_{i,j}^+\mathbf{1}_{d_{i,j}^+=0}+\alpha_{i,j}^-\mathbf{1}_{d_{i,j}^-=0})W_j(\tau)$, where $W_j$ is a standard Brownian motion.
\end{itemize}
Moreover, if for all $j=1,\dots,p$, $E(\epsilon_j(0)^{2\alpha})<\infty$ with $\alpha>1\vee (1+2d_{max})^{-1}$ where $d_{max}=max_i\{d_i\}$, then $S_n(.)$ converges towards the $p$-mfBm $X(.)$ in the Skorohod space $\mathcal D([0,1])$.
\end{prop}

\begin{proof}[Sketch of proof]
 We focus on the convergence in law of $S_n(\tau)$ to $X(\tau)$, for a fixed $\tau$ in $\R$, the finite dimensional convergence is deduced in the same way. We set for simplicity $\tau=1$.

According to the Cram\'er-Wold device, for any vector $(\lambda_1,\dots,\lambda_p)\in\R^p$,
we must show that $\lambda'S_n(1)$ converges in law to $\lambda'X(1)$. We may rewrite $\lambda'S_n(1)$ as a sum of discrete stochastic integrals (see \cite{surgailis03} and \cite{BV}) :
\begin{align}\label{sumint}
\lambda'S_n(1)&=\sum_{i=1}^p \lambda_i n^{-d_i-(1/2)}\sum_{t=1}^{n}Z_i(t)\nonumber\\
&=\sum_{i=1}^p \lambda_i \sum_{j=1}^p \sum_{k\in\Z} n^{-d_i-(1/2)} \sum_{t=1}^{n} (\psi_{i,j}^+ (t-k)+\psi_{i,j}^-(t-k)) \epsilon_j(k)\nonumber\\
&=\sum_{i=1}^p \lambda_i \sum_{j=1}^p \int_{\R}(f_{i,j,n}^+(x)+f_{i,j,n}^-(x))W_{j,n}(dx),
\end{align}
where the stochastic measures $W_{j,n}$, $j=1,\dots,p$ are defined on
finite intervals $C$ by $$W_{j,n}(C)=n^{-1/2}\sum_{k/n\in
  C}\epsilon_j(k),$$ and where $f_{i,j,n}^+$, $f_{i,j,n}^-$ are piecewise constant functions defined as follows:  denoting  $\lceil x\rceil$ the smallest integer not less than $x$, we have   for all $x\in\R$, $$f_{i,j,n}^+(x)=n^{-d_i} \sum_{t=1}^{n} \psi_{i,j}^+ (t-\lceil nx\rceil),$$ respectively $f_{i,j,n}^-(x)=n^{-d_i} \sum_{t=1}^{n} \psi_{i,j}^- (t-\lceil xn\rceil)$.\\

The following lemma states the convergence of a linear combination of discrete stochastic integrals as in (\ref{sumint}). A function is said $n$-simple if it takes a finite number a nonzero constant values on intervals $(k/n,(k+1)/n]$, $k\in\Z$.
\begin{lem}\label{cvint} Let $(f_{1,n},\cdots,f_{p,n})_{n\in\N}$ be a sequence of $p$ $n$-simple functions in $L^2(\R)$. If for any $j=1,\dots,p$, there exists $f_j\in L^2(\R)$ such that $\int_{\R}|f_{j,n}(x)-f_j(x)|^2dx\to 0$, then $\sum_{j=1}^p \int_{\R}f_{j,n}(x)W_{j,n}(dx)$ converges in law to $\sum_{j=1}^p \int_{\R}f_{j}(x)W_{j}(dx)$, where the $W_j$'s are independent standard Gaussian random measures.
\end{lem}
When $p=1$, this lemma is proved in \cite{surgailis82}. The case $p=2$ is considered in \cite{BV} and the extension to $p\geq 3$ is straightforward.

From Lemma \ref{cvint} and (\ref{sumint}), it remains to show that
$$\lim_{n\to\infty}\int_{\R}\left|f_{i,j,n}^\pm(x)-\frac{\alpha_{i,j}^\pm}{d_i} \left( (1-x)_\pm^{d_i}-
  (-x)_\pm^{d_i}\right)\mathbf{1}_{d_{i,j}^\pm=d_i}\right|^2dx = 0,$$
where we agree that $d_i^{-1}((1-x)_\pm^{d_i}- (-x)_\pm^{d_i})=\mathbf{1}_{[0,1]}(x)$ when $d_i=0$. Below, we only consider the pointwise convergence of $f_{i,j,n}^\pm(x)$, for $x\in\R$, when $d_{i,j}^\pm=d_i$. The convergence in $L^2$ is then deduced from the dominated convergence theorem (see \cite{surgailis82}, \cite{surgailis03}, \cite{BV} for details). It also follows easily that, when $d_{i,j}^\pm<d_i$, $\int_{\R}|f_{i,j,n}^\pm(x)|^2dx\to0$.

Under assumption $(i)$, note that since $d_i>0$, $(1-x)_\pm^{d_i}- (-x)_\pm^{d_i}=d_i\int_0^1 (t-x)_\pm^{d_i-1}dt$.  We have, for any $x\in\R$,
\begin{align*}f_{i,j,n}^\pm(x)&=n^{-d_i} \sum_{t=1}^{n} \psi_{i,j}^\pm (t-\lceil nx\rceil)\\
&=n^{-d_i}\int_0^n \psi_{i,j}^\pm (\lceil t \rceil-\lceil nx\rceil)dt\\
&=n^{-d_i}\int_0^n\left(\alpha_{i,j}^\pm + o(1)\right) (\lceil t \rceil-\lceil nx\rceil)_\pm^{d_{i}- 1}dt\\
&=\int_0^1\left(\alpha_{i,j}^\pm + o(1)\right) \left(\frac{\lceil nt \rceil-\lceil nx\rceil}{n}\right)_\pm^{d_{i}- 1}dt\longrightarrow \alpha_{i,j}^\pm \int_0^1 (t-x)_\pm^{d_i-1}dt.
\end{align*}

Under assumption $(ii)$, $d_i<0$. When $x\leq0$, $(1-x)_+^{d_i}- (-x)_+^{d_i}=d_i\int_0^1 (t-x)^{d_i-1}dt$ and the convergence of $f_{i,j,n}^+(x)$ can be proved as above. When $x\geq 1$, $(1-x)_+^{d_i}- (-x)_+^{d_i}=0=f_{i,j,n}^+(x)$. When $0\leq x\leq 1$, $(1-x)_+^{d_i}- (-x)_+^{d_i}=-d_i\int_1^{+\infty} (t-x)^{d_i-1}dt$ and, since $\sum_{k\in\Z} \psi_{i,j}^+(k)=0$, we have
\begin{align*}
f_{i,j,n}^+(x)&=n^{-d_i} \sum_{t=\lceil nx\rceil}^{n} \psi_{i,j}^+ (t-\lceil nx\rceil)=n^{-d_i} \sum_{t=0}^{n-\lceil nx\rceil}\psi_{i,j}^+ (t)=-n^{-d_i}\sum_{t>n-\lceil nx\rceil}\psi_{i,j}^+ (t).
\end{align*}
Therefore,
\begin{align*}
f_{i,j,n}^+(x)&=-n^{-d_i}\int_{n-\lceil
  nx\rceil}^{+\infty}\left(\alpha_{i,j}^+ + o(1)\right) (\lceil t
\rceil)^{d_{i}- 1}dt\\&=-\int_{1-\frac{\lceil
    nx\rceil}{n}}^{+\infty}\left(\alpha_{i,j}^+ +
  o(1)\right)\left(\frac{\lceil nt\rceil}{n}\right)^{d_{i}- 1}dt\\ 
& \longrightarrow
-\alpha_{i,j}^+\int_{1-x}^{+\infty}t^{d_{i}-
  1}dt=-\alpha_{i,j}^+\int_1^{+\infty} (t-x)^{d_i-1}dt 
\end{align*}
This proves $f_{i,j,n}^+(x)\to d_i^{-1}\alpha_{i,j}^+((1-x)_+^{d_i}- (-x)_+^{d_i})$, for any $x\in\R$, under assumption $(ii)$. The same scheme may be used to prove that $f_{i,j,n}^-(x)\to d_i^{-1}\alpha_{i,j}^-((1-x)_-^{d_i}- (-x)_-^{d_i})$ under assumption $(ii)$, noting that
$$(1-x)_-^{d_i}- (-x)_-^{d_i}=\begin{cases} 0 &\textrm{when}\quad x\leq0, \\ -d_i\int_{-\infty}^0(t-x)^{d_i-1}dt &\textrm{when} \quad 0\leq x\leq1, \\d_i\int_0^1 (t-x)_-^{d_i-1}dt  &\textrm{when} \quad x>1.\end{cases}$$

Under assumption $(iii)$, $$f_{i,j,n}^\pm(x)=\sum_{t=1}^n\psi_{i,j}^\pm (t-\lceil nx\rceil)=\sum_{t=1-\lceil nx\rceil}^{n-\lceil nx\rceil}\psi_{i,j}^\pm (t).$$ Since $\sum_{t\in\Z}\psi_{i,j}^\pm (t)<\infty$, $f_{i,j,n}^\pm(x)\to0$ for all $x\notin[0,1]$. When $x\in[0,1]$, we have $f_{i,j,n}^\pm(x)\to \alpha_{i,j}^\pm$.

Therefore, the first claim of the theorem is proved, i.e. the convergence in law of the finite dimensional distribution of $(S_n(\tau))_{\tau\in\R}$ to $(X(\tau))_{\tau\in\R}$. To extend this convergence to a functional convergence in $\mathcal D([0,1])$, it remains to show tightness of the sequence $(S_n(\tau))_{\tau\in[0,1]}$. This follows exactly from the same arguments as in the proof of Theorem 1.2 in \cite{BV}.
\end{proof}
\section{Synthesis of the mfBm}\label{simulation:sec}

\subsection{Introduction}

The exact simulation of the fractional Brownian motion has been a
question of great interest in the nineties. This may be done by
generating a sample path of a fractional Gaussian noise. An important
step towards efficient simulation was obtained after the work of Wood
and Chan \cite{WoodC94} about the simulation of arbitrary
stationary Gaussian sequences with prescribed covariance function. The
technique relies upon the embedding of the covariance matrix into a circulant matrix, a square root of which is easily calculated using the discrete Fourier transform. This leads to a very efficient algorithm, both in terms of computation time and storage needs.
Wood and Chan method is an exact simulation method provided that
the circulant matrix is semidefinite positive, a property that is not
always satisfied. However, for the fractional Gaussian noise, it can
be proved that the circulant matrix is definite positive for all $H
\in (0,1)$, see \cite{Craig03,DietN97}.

In \cite{ChanW99}, Wood and Chan extended their method and provided a more general algorithm adapted to multivariate stationary  Gaussian processes. The main characteristic of this method is that if a certain condition for a familiy of Hermitian matrices holds then the algorithm is exact in principle, i.e. the simulated data have the true covariance. We present hereafter the main ideas, briefly describe the algorithm and propose some examples.

\begin{rem}
 Other approaches could have been undertaken (see
\cite{bardet} for a review in the case $p=1$). Approximate simulations can be done by discretizing the moving-average or spectral stochastic integrals~(\ref{eq:1}) or~(\ref{Xma}). \cite{Cham95} also proposed an approximate method based on the spectral density matrix of the increments for synthesizing multivariate Gaussian time series. Thanks to Proposition~\ref{prop-spd}, this could also be envisaged for the mfBm.
\end{rem}

\subsection{Method and algorithm}

 For two arbitrary matrices $A=(A_{j,k})$ and $B$, we use $A\otimes B$ to denote the Kronecker product of $A$ and $B$ that is the block matrix $(A_{j,k} B)$.

Let $\Delta X := \Delta_1 X$ denotes the increments of size 1
($\delta=1$) of a mfBm. We have $\Delta X = (\Delta X(t))_{t\in\R}=
((\Delta X_1(t),\ldots,\Delta X_p(t))')_{t\in\R}$.
The aim is to simulate a realization of a multivariate fractional
Gaussian noise discretized at times $j=1,\ldots,n$, that is a
realization of $(\Delta X(1),\ldots,\Delta X(n))$. Then a realization
of the discretized mfBm will be easily obtained.

We denote by
$\Delta X^{(n)}$ the merged vector $\Delta X^{(n)} =(\Delta
X(1)',\ldots,\Delta X(n)')'$ and by $\mathbb{G}$ its covariance
matrix. $\mathbb{G}$ is the $np\times np$ Toeplitz block matrix $\mathbb{G} =
 (G(|i-j|))_{i,j=1,\dots,n}$ where for $h=0,\ldots,n-1$, $G(h)$ is the
 $p\times p$ matrix given by $G(h):=\left( \gamma_{j,k}(h)
 \right)_{j,k=1,\ldots,p}$. The simulation problem can be viewed as
 the generation of a random vector following a
 $\mathcal{N}_{np}(0,\mathbb G)$. This may be done by computing $\mathbb G^{1/2}$ but the complexity of such a procedure is $\mathcal{O}(p^3n^2)$ for block Toeplitz matrices.
To overcome this numerical cost, the idea is to embed $\mathbb G$ into the block circulant matrix $C=\mathrm{circ} \{ C(j), j=0,\ldots,m-1\}$, where $m$ is a power of 2 greater than $2(n-1)$ and where each $C(j)$ is the $p\times p$ matrix defined by
\begin{equation} \label{eq-C}
C(j)= \left\{ \begin{array}{ll}
G(j) & \mbox{ if } 0\leq j < m/2 \\
\frac 12 \big(G(j)+{G(j)}^\prime\big) & \mbox{ if } j=m/2 \\
G(j-m) & \mbox{ if } m/2 < j \leq m-1.
\end{array} \right.
\end{equation}
Such a definition ensures that $C$ is a symmetric matrix with nested
block circulant structure and that $\mathbb G=\{ C(j),
j=0,\ldots,n-1\}$ is a submatrix of $C$. Therefore, the simulation of
a $\mathcal{N}_{np}(0,\mathbb G)$ may be achieved by taking the $n$ ``first'' components of a vector $\mathcal{N}_{mp}(0,C)$, which is done by computing $C^{1/2}$. The last problem is more simple since one may exploit the circulant characteristic of $C$: there exist $m$ Hermitian matrices $B(j)$ of size $p\times p$ such that the following decomposition holds
\begin{equation} \label{eq-WC}
C = (J \otimes I_p) \; \mathrm{diag}(B(j),j=0,\ldots,m-1) \; (J^* \otimes I_p),
\end{equation}
where $Q$ is the $m\times m$ unitary matrix defined for $j,k=0,m-1$ by $J_{j,k}=e^{-2\i\pi jk/m}$. The computation of $C^{1/2}$ is much less expensive than the computation of $\mathbb G^{1/2}$ since, as in the one-dimensional case ($p=1$), (\ref{eq-WC}) will allow us to make use of the Fast Fourier Transform (FFT) which considerable reduces the complexity.

Now, the algorithm proposed by Wood and Chan may be described through the following steps. Let $m$ be a power of 2 greater than $2(n-1)$.

\noindent\underline{Step 1.} For $1\leq u\leq v\leq p$ calculate for $k=0,\ldots, m-1$
$$
B_{u,v}(k) = \sum_{j=0}^{m-1} C_{u,v}(j) e^{-2
\i\pi jk/m}
$$
where $C_{u,v}(j)$ if the element $(u,v)$ of the matrix $C(j)$ defined by~(\ref{eq-C}) and set $B_{v,u}(k)=B_{u,v}(k)^*$.\\

\noindent\underline{Step 2.} For each $j=0,\ldots,m-1$ determine a unitary matrix $R(j)$ and real numbers $\xi_u(j)$ ($u=1,\ldots,p$) such that ${B}(j)=R(j) \; \mbox{diag}(\xi_1(j),\ldots,\xi_p(j)) \; R(j)^*$.\\

\noindent\underline{Step 3.} Assume that the eigenvalues $\xi_1(j),\ldots,\xi_p(j)$ are non-negative (see Remark~\ref{rem-condition}) and define $\widetilde{B}(j)=R(j) \; \mbox{diag}(\sqrt{\xi_1(j)},\ldots,\sqrt{\xi_p(j)}) \; R(j)^*$. \\

\noindent\underline{Step 4.} For $j=0,\ldots,m/2$ generate independent vectors $U(j),V(j) \sim \mathcal{N}_p(0,I)$ and define
$$
Z(j)= \frac{1}{\sqrt{2m}}\times\left\{ \begin{array}{ll}
\sqrt{2} U(j) & \mbox{for } j=0,\frac m2 \\
U(j)+ \i V(j) & \mbox{for } j=1,\ldots,\frac m2-1 ,\\
\end{array}\right.
$$
let $Z(m-j)=\overline{Z}(j)$ for $j=\frac m2+1,\ldots,m-1$ and set $W(j):=\widetilde{B}(j) Z(j)$.\\

\noindent\underline{Step 5.} For $u=1,\ldots,p$ calculate for $k=0,\ldots,m-1$
$$
\Delta X_u(k) = \sum_{j=0}^{m-1} W_u(j) e^{-2\i\pi jk/m}
$$
and return $\big\{ \Delta X_u(k), 1\leq u \leq p, k=0,\ldots,n-1\big\}$. \\

\noindent\underline{Step 6.} For $u=1,\ldots,p$ take the cumulative sums $\Delta X_u$ to get the $u-th$ component $X_u$ of a sample path of a mfBm.

\begin{figure}[ht]
  \begin{center}
    \begin{tabular}{lll}
      {\includegraphics[scale=.275]{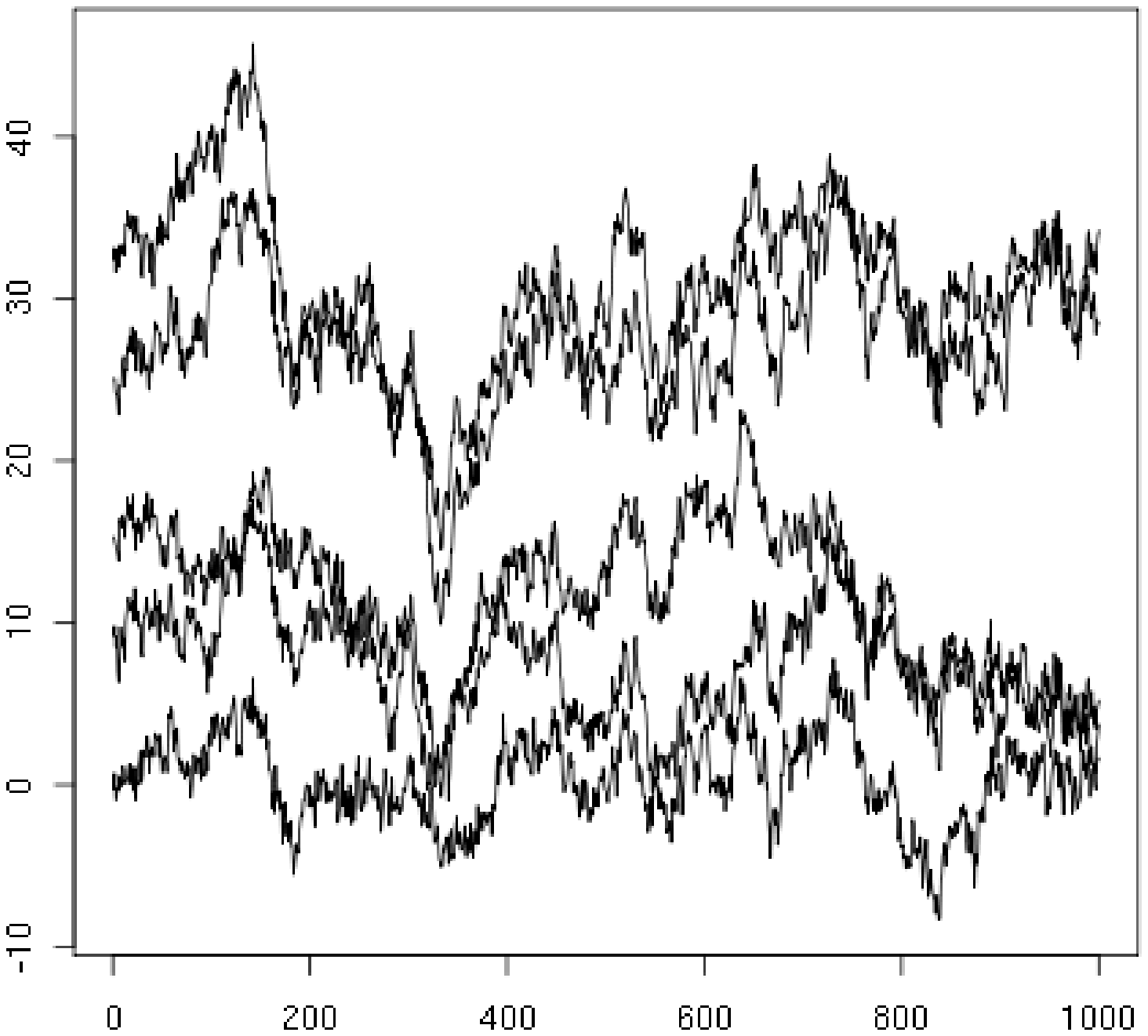}} & {\includegraphics[scale=.275]{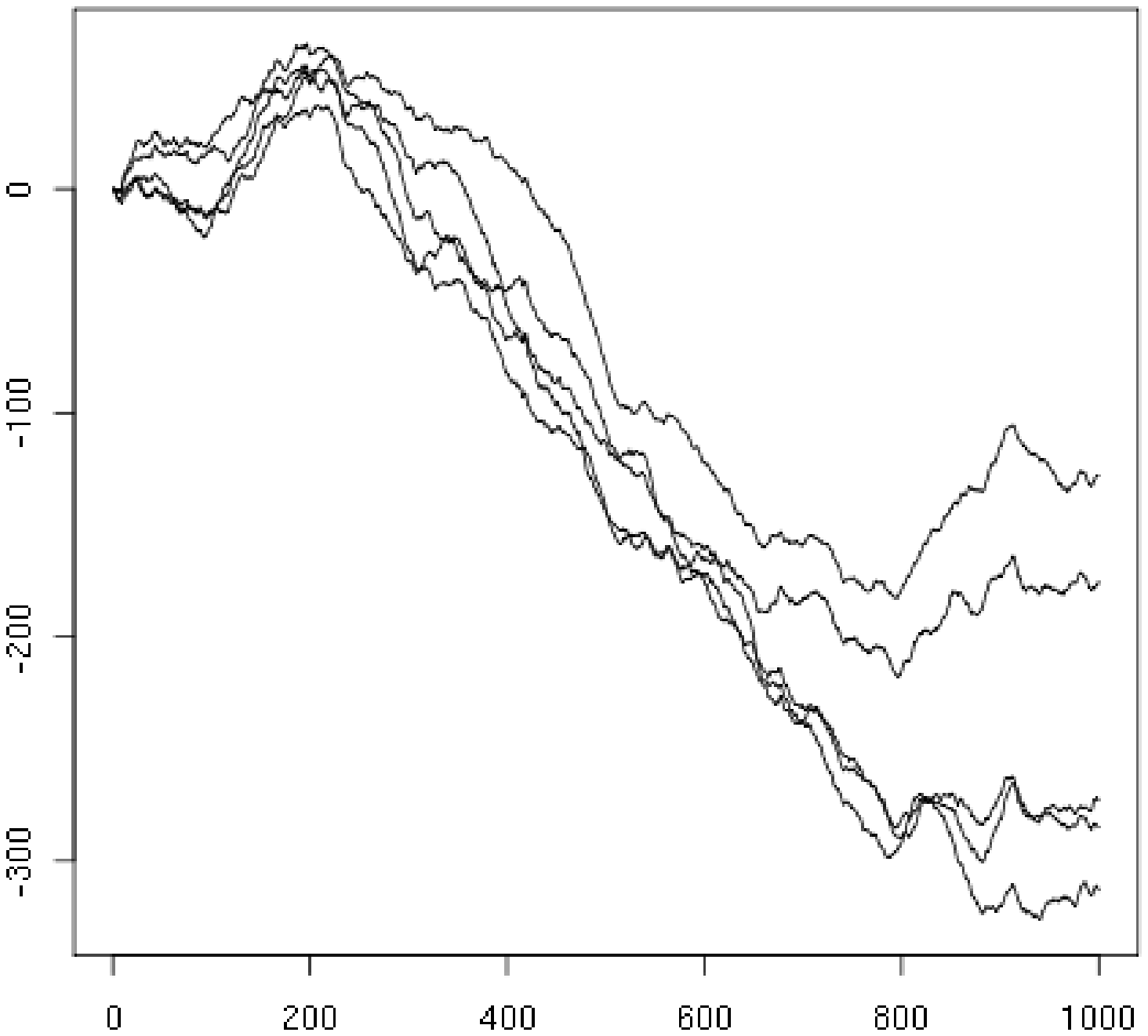}} & {\includegraphics[scale=.275]{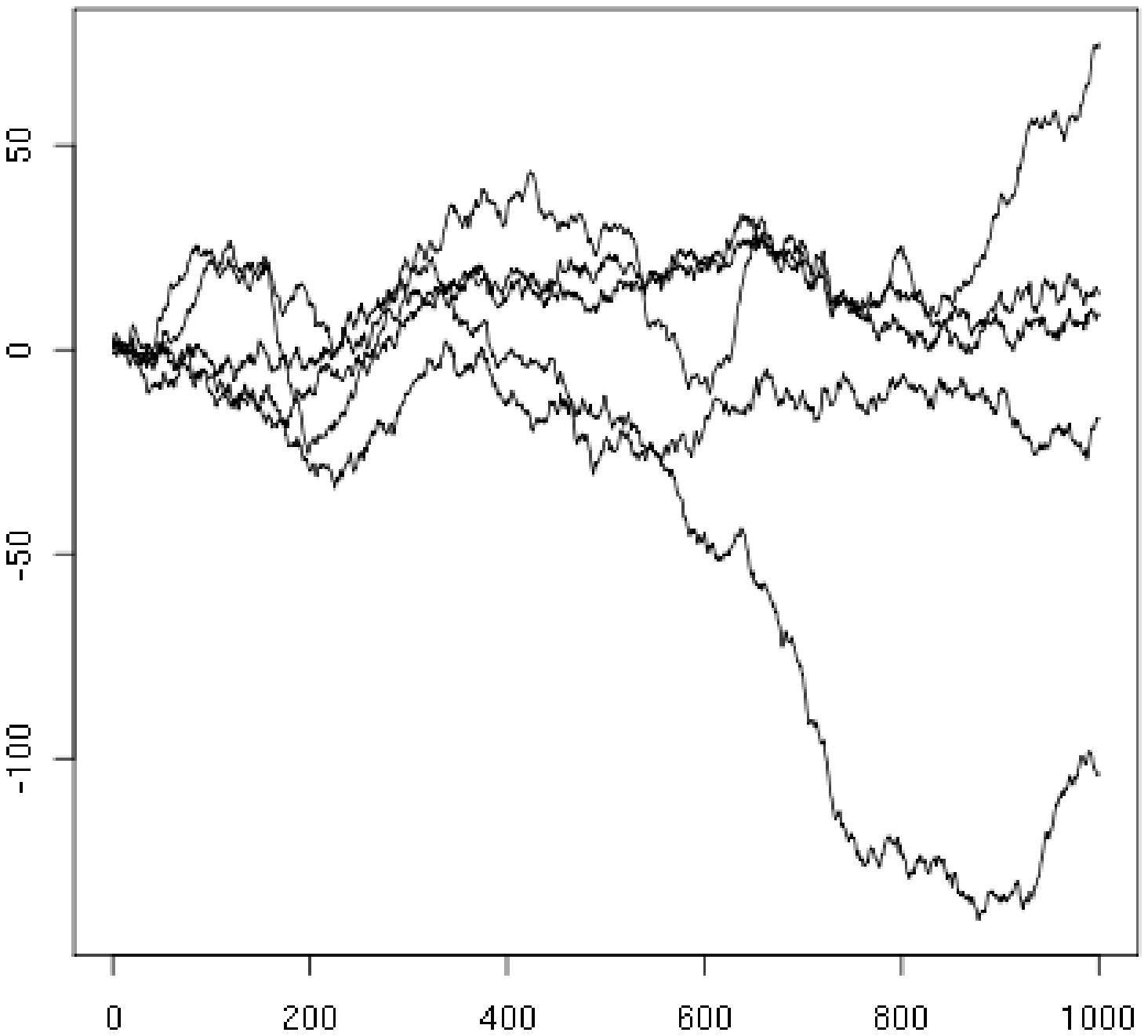}} \\
      {\includegraphics[scale=.275]{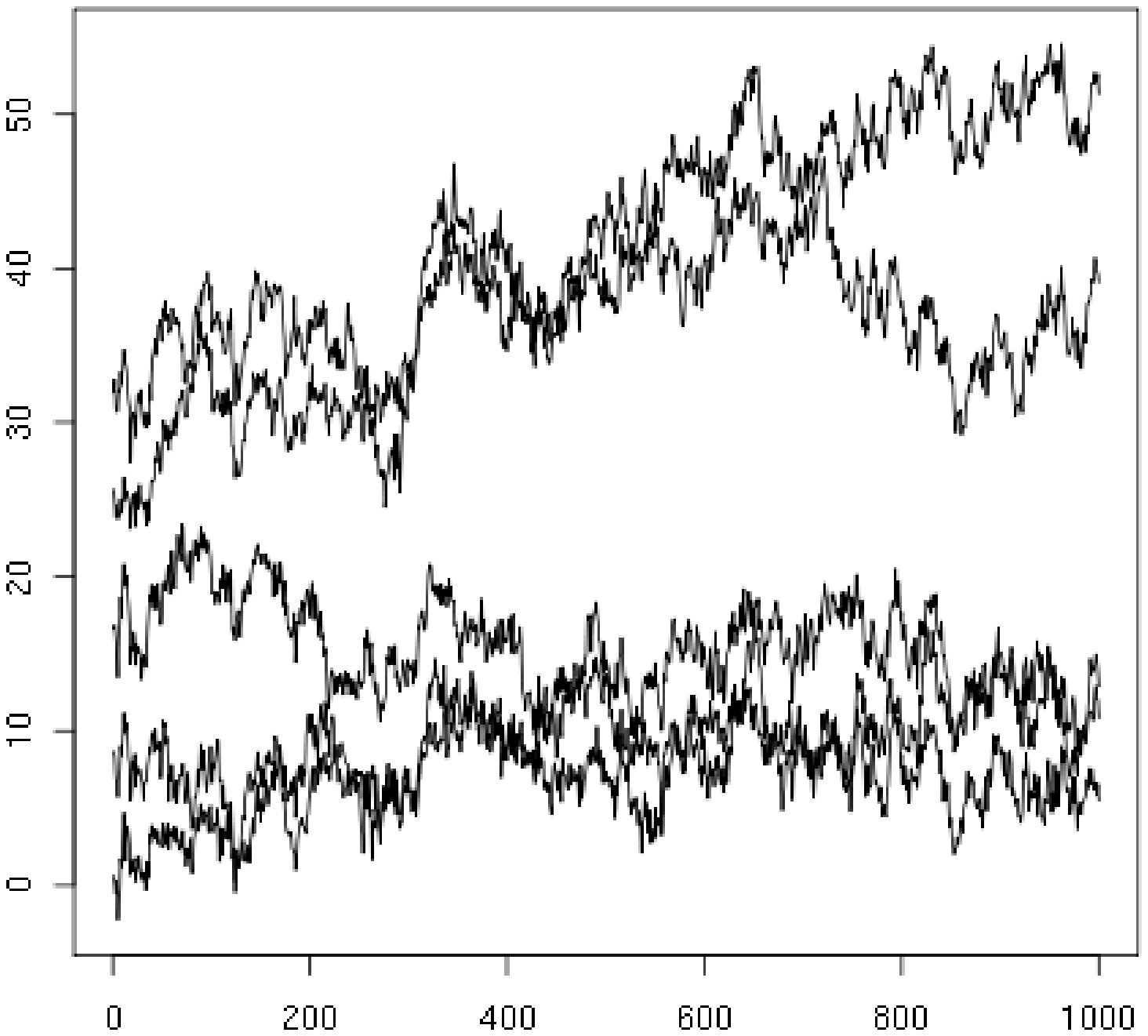}} & {\includegraphics[scale=.275]{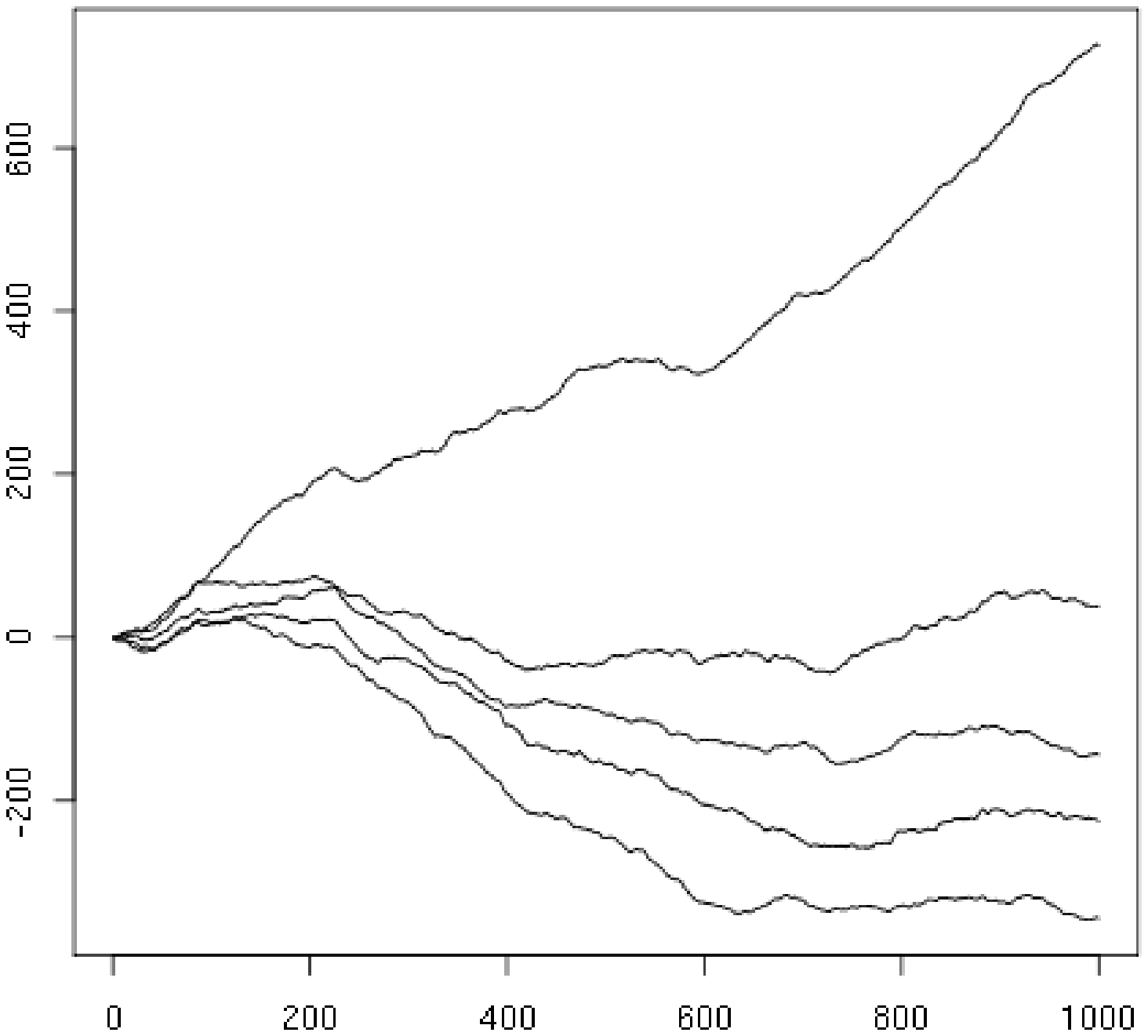}} & {\includegraphics[scale=.275]{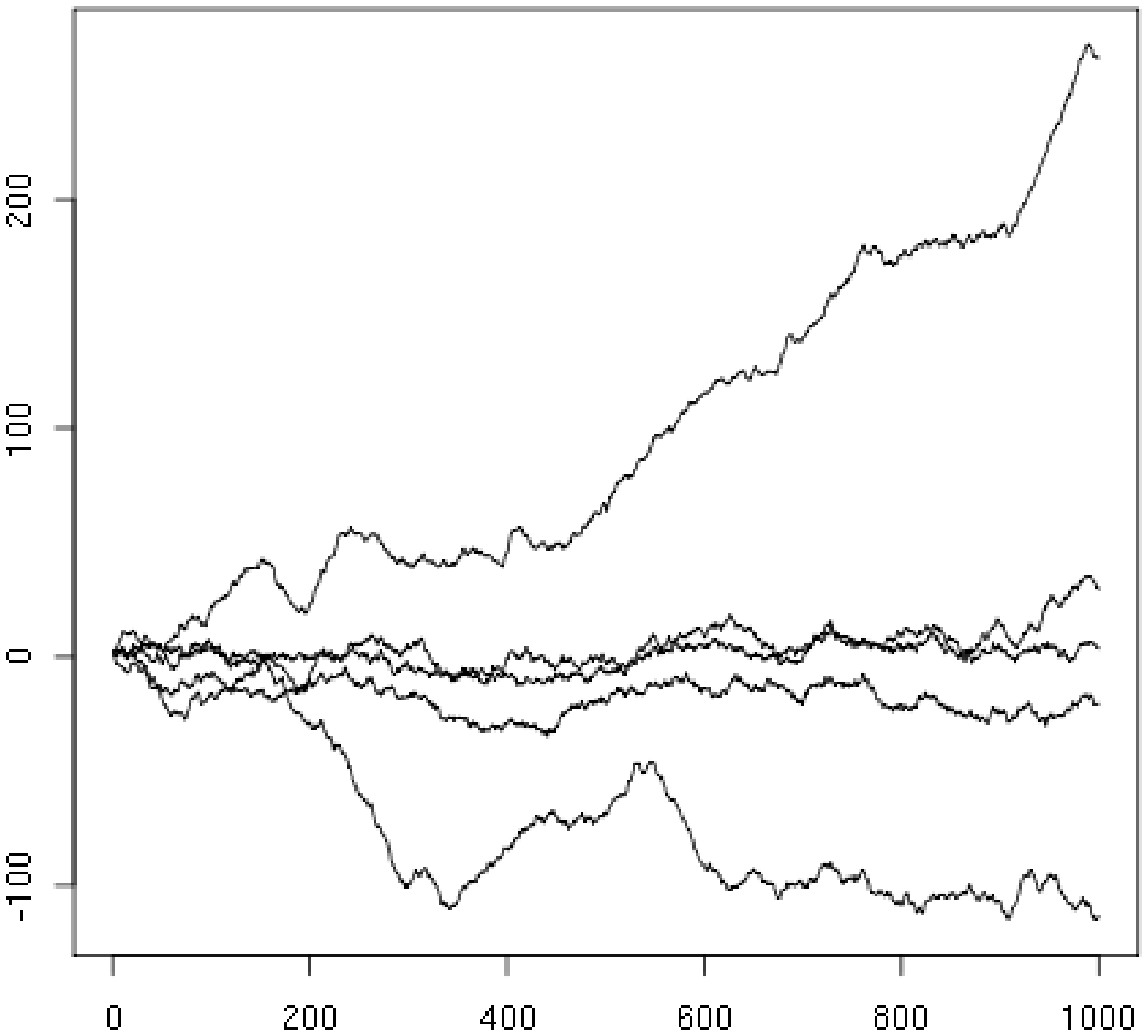}} \\
      {\includegraphics[scale=.275]{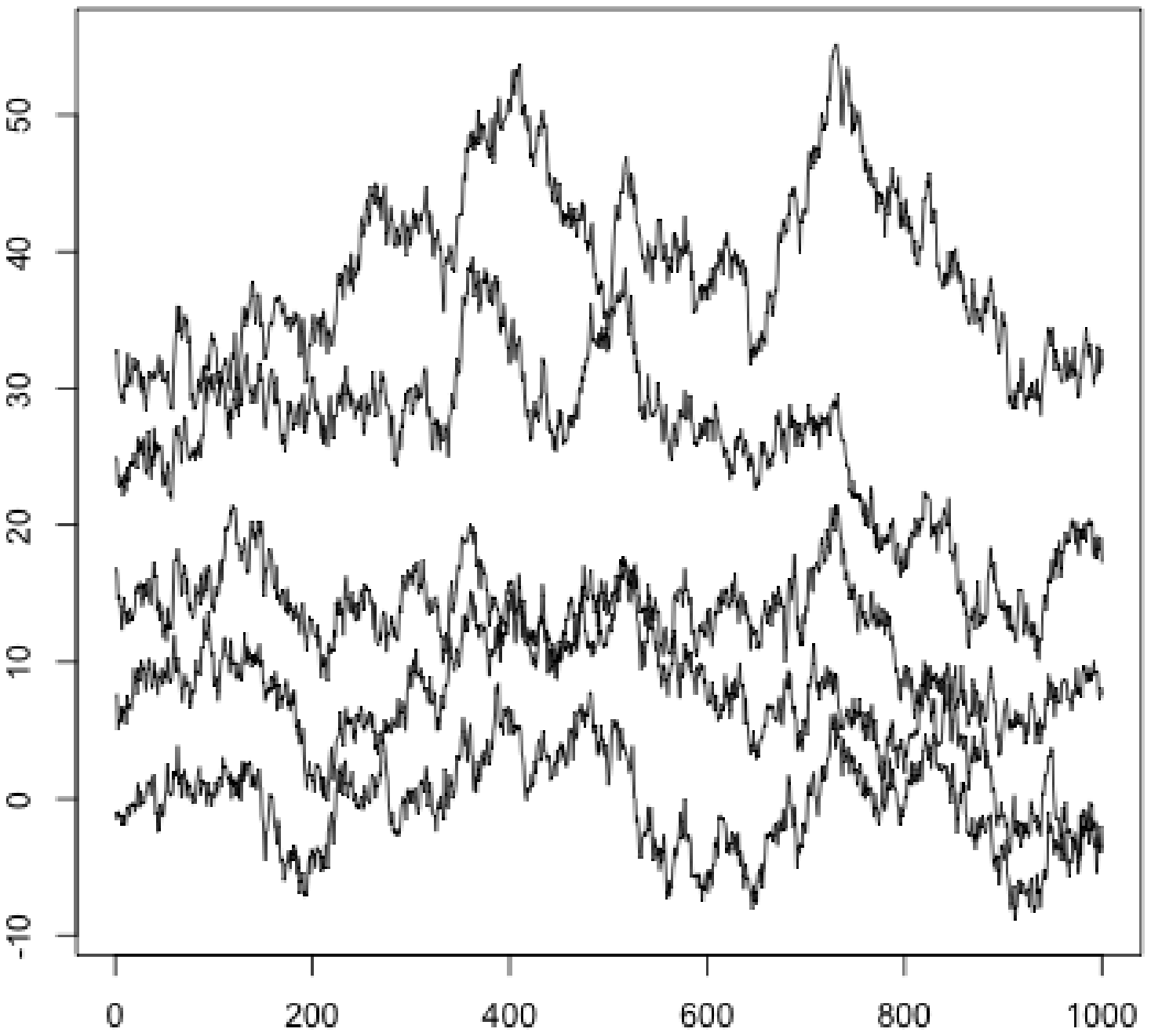}} &
      {\includegraphics[scale=.275]{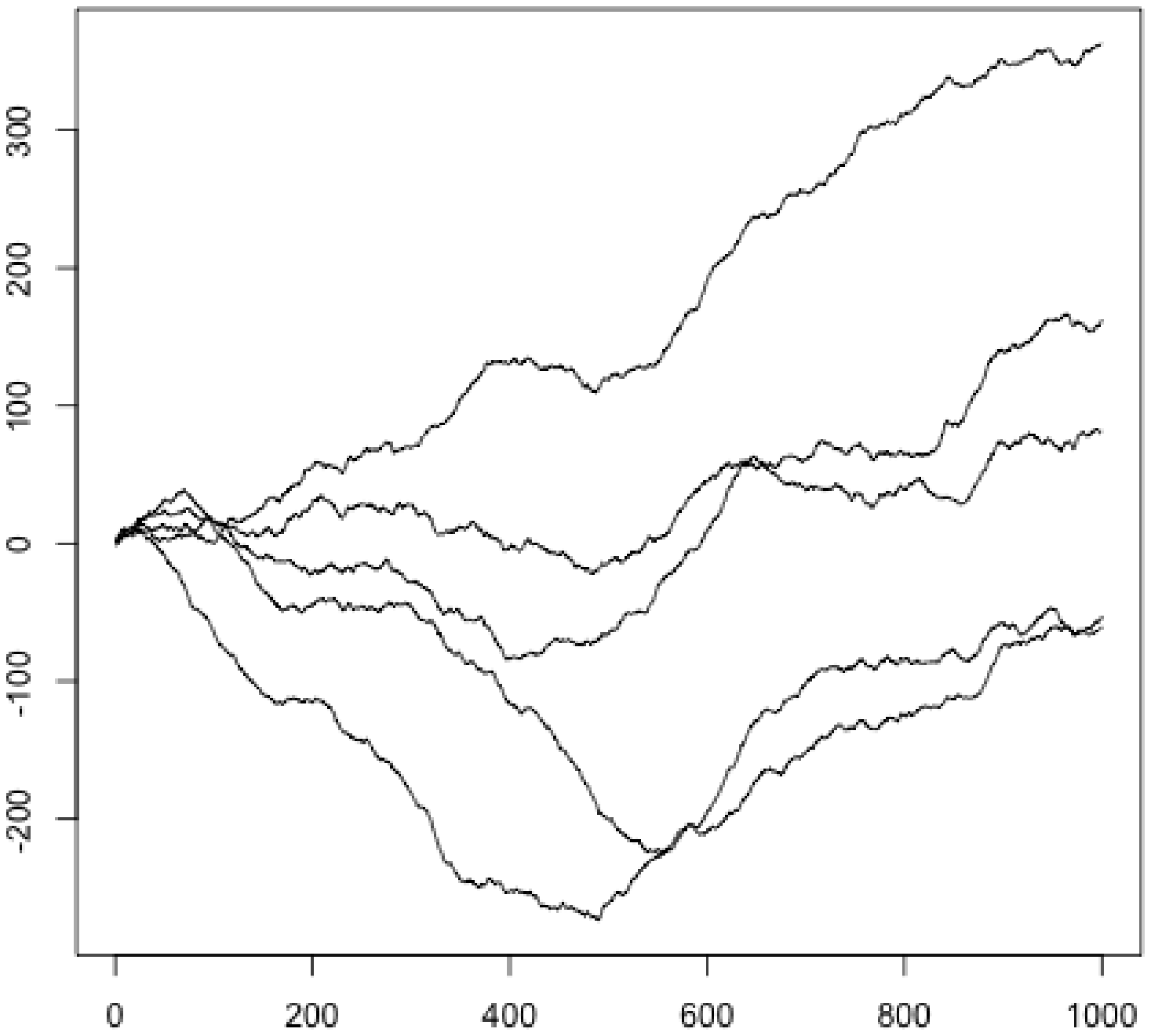}} &
      {\includegraphics[scale=.275]{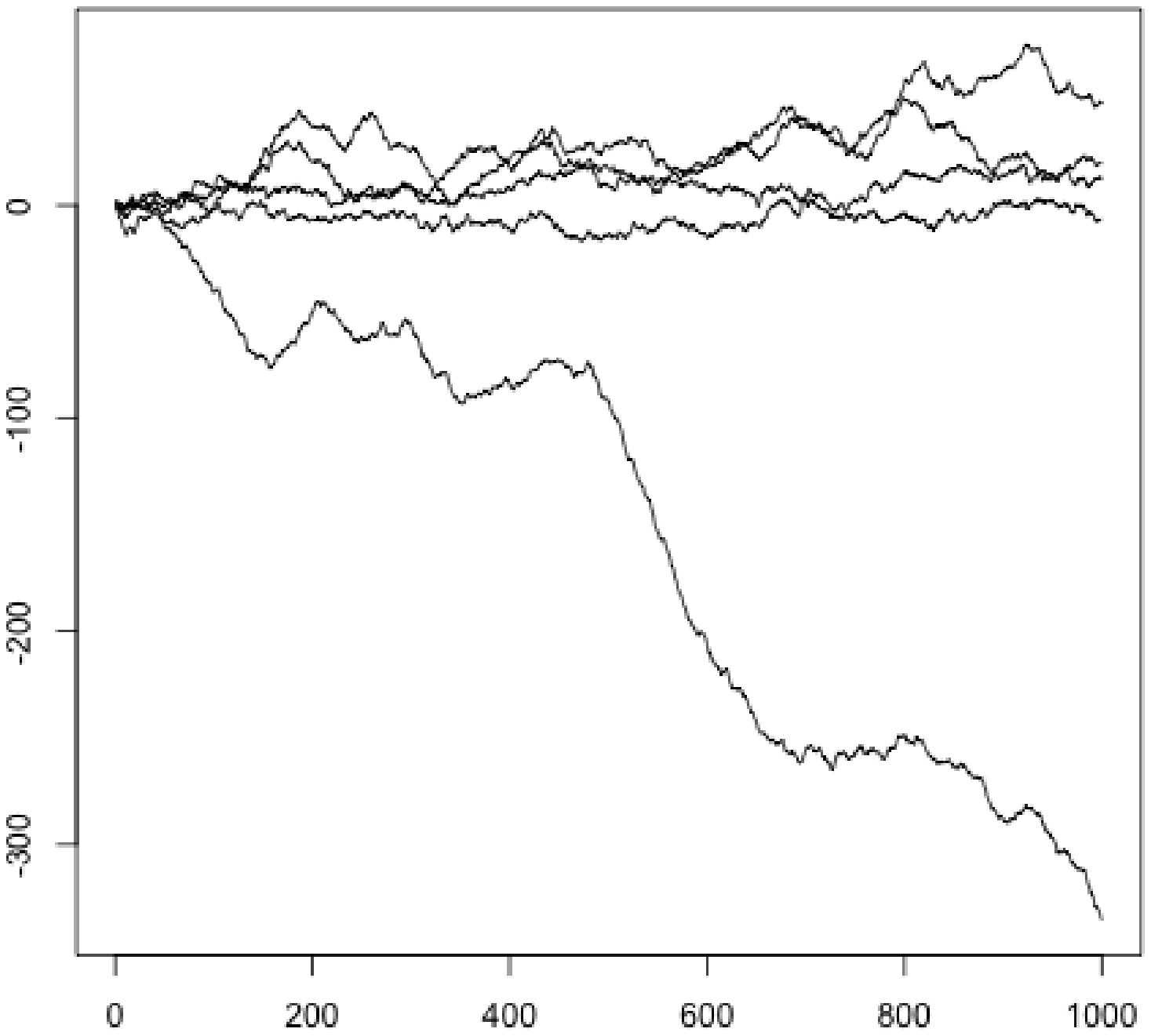}}
    \end{tabular}
  \end{center}
\caption{
Examples of sample paths of length $n=1000$ normalized causal
 (top), well-balanced (middle) and general (bottom) multivariate
 fractional Brownian motion with $p=5$ components. The Hurst
 parameters are equally spaced in $[0.3,0.4]$ (left), $[0.8,0.9]$
 (middle) and $[0.4,0.8]$ (right). The correlation parameters
 $\rho_{i,j}$ are all set to $0.4$ (left, middle) and $0.3$
 (right). For the general mfBm (bottom), the parameters $\eta_{i,j}$
 are set to {$0.1/ (1-H_i-H_j)$} and the parameters $\widetilde{\eta}_{i,j}$ to $0.1$. Note that the existence condition discussed in Proposition~\ref{prop-exist} and the condition in Step~3 of the algorithm are satisfied for these different choices of parameters. For convenience, the sample paths of the left column have been artificially shifted for visibility.}
\label{fig-mfbm}
\end{figure}

Figure~\ref{fig-mfbm} gives some examples of sample paths of mfBm's
simulated with this algorithm.

\begin{rem} Let us discuss the computation cost of the most expensive steps, that is steps 1, 2 and 5. Step 1 requires $\frac{p(p+1)}2$ applications of the FFT of signals of length $m$, Step 2 needs $m$ diagonalisations of $p\times p$ Hermitian matrices and Step 5 requires $p$ applications of the FFT of signals of length $m$. Therefore, the total cost, $\kappa(m,p)$ equals
$$\kappa(m,p) = \mathcal{O}\left( \frac{p(p+1)}2 m\log m\right) + \mathcal{O}(m p^3) + \mathcal{O}(p m\log m).$$
\end{rem}

\begin{rem} \label{rem-condition}
The crucial point of the previous algorithm lies in the non-negativity
of the eigenvalues $\xi_1(j),\ldots, \xi_p(j)$ for any
$j=0,\ldots,m-1$. In the one-dimensional case (when $p=1$) Steps 2 and
3 disappear, and in Step 1, $B_{11}(k)$ corresponds to the $k-$th
eigenvalue of the circulant matrix $C_{11}$ with first line defined by
$C_{11}(j)=\gamma_{11}(j)$ for $0\leq j \leq m/2$ and $\gamma(m-j)$
for $j=m/2+1,\ldots,m-1$. For the fractional Gaussian noise, it has
been proved by Craigmile \cite{Craig03} for $H<1/2$, and by Dietrich
and Newsam \cite{DietN97} for $H>1/2$ that such a matrix is semidefinite-positive for any $m$ (and so for the first power of 2 greater than $2(n-1)$). In the more general case $p>1$, the problem is much more complex: the quantities $B_{u,v}(k)$ are not necessarily real, and the establishment of a condition of positivity for the matrix $B_{u,v}(k)$ does not seem obvious. When the condition in Step 3 does not hold, Wood and Chan suggest to either increase the value of $m$ and restart Steps 1,2 or to truncate the negative eigenvalues to zero which leads to an approximate procedure. These problems are not addressed in this paper.
Let us assert that for the simulation examples presented in
Figure~\ref{fig-mfbm}, we have observed that this condition is
satisfied for $m$ equal to the first power of 2.
\end{rem}

\bibliographystyle{apalike}

\end{document}